\documentclass[sts]{imsart}
\usepackage{amsmath}
\usepackage{amsfonts}
\usepackage{amssymb}
\usepackage{appendix}
\usepackage{amsthm}
\usepackage{footnote}
\usepackage{hyperref}
\usepackage{fullpage}
\usepackage{graphicx}
\usepackage{epstopdf}
\usepackage[numbers,colon]{natbib}

\usepackage{enumerate}
\newtheorem{Thm}{Theorem}[section]
\newtheorem{Cla}[Thm]{Claim}

\newtheorem{Lem}[Thm]{Lemma}

\makeatletter
\def\blfootnote{\xdef\@thefnmark{}\@footnotetext}
\makeatother

\theoremstyle{definition}
\newtheorem{Def}[Thm]{Definition}

\theoremstyle{remark}
\newtheorem{Rem}[Thm]{\bf{Remark}}

\newcommand{\ConvD}{\overset{d}{\rightarrow}}
\newcommand{\ConvP}{\overset{p}{\rightarrow}}
\newcommand{\ConvFDD}{\overset{f.d.d.}{\longrightarrow}}

\newcommand{\Cov}{\mathrm{Cov}}
\newcommand{\Var}{\mathrm{Var}}
\newcommand{\E}{\mathbb{E}}
\newcommand{\mathbd}{\boldsymbol}

\begin{document}

\begin{frontmatter}
\title{How the instability of ranks under long memory   affects large-sample inference}
\runtitle{Instability of ranks under long memory}

\begin{aug}
\author{\fnms{Shuyang} \snm{Bai} \ead[label=e1] {bsy9142@uga.edu}},
\author{\fnms{Murad S.} \snm{Taqqu}\thanksref{t1}\ead[label=e2]{murad@bu.edu}}

\thankstext{t1}{Supported by NSF grant  DMS-1309009  at Boston University}
\runauthor{S. Bai and M.S. Taqqu}

\affiliation{University of Georgia and Boston University}

\address{Shuyang Bai: Department of Statistics,  University of Georgia, Athens, GA 30602, US \printead{e1}.}

\address{Murad S. Taqqu: Department of Mathematics and Statistics,
Boston University,
Boston, MA 02215, US \printead{e2}.}

\end{aug}

\begin{abstract}
Under long memory,  the limit theorems for normalized sums of random variables   typically involve a  positive integer called   ``Hermite rank''.  There is a different limit for each Hermite rank.    From a statistical point of view, however, we argue  that  a rank other than one is unstable, whereas, a  rank equal to one is stable.  We provide empirical evidence supporting this argument. This has important consequences. Assuming  a  higher-order rank when it is not really there usually results in  underestimating  the  order of the fluctuations of the statistic of interest. We illustrate this through various examples involving the sample variance, the empirical processes and the Whittle estimator.
\end{abstract}


\end{frontmatter}

\blfootnote{
\begin{flushleft}
\textbf{2010 AMS Classification:} 62M10, 	60F05  
\end{flushleft}
}
\section{Introduction}
Suppose that $\mathcal{D}$ is a data set, and one has a statistical model for $\mathcal{D}$ which involves a random stationary sequence $\{X(n)\}$, referred to as noise. 
Let $T=T(\mathcal{D})$ be a sample statistic of interest.
Deriving the asymptotic distribution for the statistic $T$ as the sample size tends to  infinity is a standard practice in large sample inference. The asymptotic distribution is useful for reporting confidence intervals, conducting hypothesis tests, etc.

When the construction of the statistic $T$ involves summing the data  and if the stationary noise  $\{X(n)\}$ is weakly dependent, then the asymptotic distribution of $T$ is typically Gaussian in view of the Central Limit Theorem. This asymptotic distribution can also be  a functional of a Gaussian process.

The situation, however, is much more intricate when the strength of dependence in the noise increases significantly. This  strong-dependence regime, often called  \emph{long memory} or \emph{long-range dependence}, is typically characterized by the following behavior of the variance of partial sums:
\begin{equation}\label{eq:var allan}
\Var\Big[\sum_{n=1}^N X(n)\Big]\approx {N^{2H}}, \text{ as }N\rightarrow\infty.
\end{equation}
where $\approx$ means asymptotic equivalence up to some positive constant, the parameter $H\in(1/2,1)$ is called the \emph{Hurst index}\footnote{The term ``Hurst index'' is also frequently used for the self-similarity parameter of self-similar processes arising from the normalized limit of the sum of $X(n)$ (see \citet{pipiras:taqqu:2017:long}). It is also common to introduce the so-called memory parameter $d$ in the context of (\ref{eq:cov decay}),  re-expressed as $\Cov[X(n),X(0)]\approx n^{2d-1}$ ($d=H-1/2$). We  will use  $H$ throughout   in order not to switch between parameters and thus to avoid   confusion.}.
Normally when the dependence is weak, one expects $H=1/2$ in (\ref{eq:var allan}), that is, the growth of the variance is linear. The superlinear growth in (\ref{eq:var allan}) is typically due to the slow  decay of the covariance of $\{X(n)\}$:
\begin{equation}\label{eq:cov decay}
\Cov[X(n),X(0)]\approx n^{2H-2}, \text{ as }n\rightarrow\infty,
\end{equation}
where $-1<2H-2<0$.  In fact, (\ref{eq:cov decay}) is also a common characterization of long memory.
We refer the reader to the recent monographs \citet{beran:2013:long}, \citet{giraitis:koul:surgailis:2009:large}, \citet{samorodnitsky:2016:stochastic} and \citet{pipiras:taqqu:2017:long} for   comprehensive introductions to the notion long memory.

In view of (\ref{eq:var allan}), when deriving the asymptotic distribution of the sum, one needs to associate the stronger normalization $N^{-H}$ to $\sum_{n=1}^N X(n)$ rather than the standard  $N^{-1/2}$ normalization.
These   limit theorems have been applied  in many statistical studies. See Section \ref{Sec:examples} below.

This paper makes the following basic argument: while these  limit theorems  are definitely of probabilistic interest,   their immediate application  to statistical inference can lead to problems. This is because these limit  theorems can be unstable, that is, they often cease to hold when  $\{X(n)\}$ is slightly perturbed.
In particular, the  limit theorems under long memory often depend critically on an integer quantity called \emph{rank}, e.g,   the \emph{Hermite rank} in the Gaussian context. We will show that the rank  is unstable when it takes value greater than one, and it easily collapses to rank one  when there is a slight perturbation.

The notion of rank, however, is not relevant when the data is weakly dependent. We indicate that under weak dependence, limit theorems are robust against, for example, a transformation of the data. We illustrate this by considering various types of weak dependence, such as strong mixing, Gaussian subordination and Bernoulli shifts.

The paper is organized as follows.
The rank instability issue is discussed  in Section \ref{Sec:Instability}. In Section \ref{Sec:examples}, we provide some examples on how the instability of rank can affect statistical results. In Section \ref{Sec:empirical}, we carry out an empirical  study which supports the instability argument.
In contrast, we show in Section \ref{Sec:clt stable} that the central limit theorems under weak dependence are not subject to such instability issues.
   Section \ref{Sec:conclusion} contains   conclusions and suggestions. Some technical extensions  are found  in the  Appendix.

\section{The instability of  ranks  under long memory}\label{Sec:Instability}
In this section, we   introduce   the notion Hermite rank and point out its instability.  We  focus on   the simple scenario of instantaneous transformation of a Gaussian stationary process.   (The   non-instantaneous case is somewhat  technical and   is deferred to  Appendix \ref{Sec:non-inst trans}.) We then address the case where   the model involves a non-Gaussian linear (moving-average) process, where the corresponding notion of Hermite rank is called the \emph{Appell rank} or the \emph{power rank}.

Throughout the paper, the notation $a(n) \approx b(n)$ means  $\lim_{n\rightarrow\infty}a(n)/b(n)= c$ for some generic constant $0<c <\infty$ that can change from expression to expression. We note that in many places, one can  include a slowly varying function in the asymptotic relation, for example,  a logarithmic function (see, e.g., \citet{bingham:goldie:teugels:1989:regular}), but for simplicity we do not do that.

\smallskip
\subsection{Transformation of Gaussian processes}\label{Sec:inst trans}
We want to consider possibly nonlinear finite-variance transformations of Gaussian random variables. To do so,
let $Z$ be a standard normal random variable, $\gamma(dx)$ be the standard Gaussian measure $(2\pi)^{-1/2} e^{-x^2/2} dx$ on $\mathbb{R}$, and let
\[
L^2(\gamma)=\{G(\cdot):~\E G(Z)^2<\infty\}.
\]

It is well-known (see, e.g., see \citet{pipiras:taqqu:2017:long}, Proposition 5.1.3) that $\{\frac{1}{\sqrt{m!}}H_m(\cdot), m\ge 0\}$ forms an orthonormal basis of $L^2(\gamma)$, where
 $\{H_m(\cdot), ~m\ge 0\}$ are Hermite polynomials defined as $H_0(x)=1$ and
\begin{equation}\label{eq:Herm Poly}
H_m(x)=(-1)^m e^{x^2/2}\frac{d}{dx^m}e^{-x^2/2} \quad\text{ for }\quad m\ge 1.
\end{equation}
Thus $H_1(x)=x$, $H_2(x)=x^2-1$ and $H_3(x)=x^3-3x$, etc.   We can now define the Hermite rank of a function $G\in L^2(\gamma)$.

\begin{Def}\label{Def:Herm rank}
Suppose that  $G(\cdot)\in L^2(\gamma)$. Let $Z$ be a standard Gaussian random variable. The  \emph{Hermite  rank}  $k$ of $G(\cdot)$ is defined as
\begin{equation}\label{eq:def herm rank}
k=\inf\left\{m\ge 1: ~\E G(Z) H_m(Z)=\int_{\mathbb{R}}G(x)H_m(x)\gamma(dx)\neq 0 \right\}.
\end{equation}
where $H_m(\cdot)$ is the $m$-th order Hermite polynomial.
\end{Def}

\begin{Rem}
An alternative  way of defining the Hermite rank $k$ is through the starting index of the Hermite expansion of $G(\cdot)-\E G(Z)$, namely,
\begin{equation}\label{eq:def herm rank expan}
G(\cdot)-\E G(Z)=\sum_{m=k}^\infty c_m H_m(\cdot), \quad k\ge 1,
\end{equation}
for some sequence $c_m$ satisfying $c_k\neq 0$, where the series converges in the $L^2(\gamma)$-sense.  By the orthonormality of $\{\frac{1}{\sqrt{m!}}H_m(\cdot) \}$, we have
\begin{equation}\label{eq:Herm Coeff}
c_m={\E[G(Z)H_m(Z)]\over m!} ~~\text{ for }~~m\ge 0.
\end{equation}
Note that $c_0=\E G(Z)$ since $H_0(Z)=1$. Furthermore, since the Hermite polynomials $\{H_m(\cdot), 0\le m\le k\}$ form a basis for
polynomials of degree less than equal to $k$, the definition (\ref{eq:def herm rank}) can be re-expressed as
\begin{equation}\label{eq:alter herm rank def}
k=\inf\left\{m\ge 1: ~\E \Big[\big(G(Z)-\E G(Z)\big) Z^m\Big]\neq 0 \right\}.
\end{equation}
\end{Rem}

\begin{Rem}
The Hermite rank of $G(x)$ is the same as that of $G(x)+a$, for any $a\in\mathbb{R}$, since relation (\ref{eq:def herm rank expan}) involves centering.
\end{Rem}

Now suppose that  $\{X(n)\}$ is a long-memory stationary Gaussian process satisfying (\ref{eq:cov decay}). We may assume    without loss of generality that it is  standardized, that is,  $\E X(n)=0$ and $\Var[X(n)]=1$
The following lemma explains the role that the Hermite rank plays in determining the  asymptotic behavior of the covariance of the transformed sequence $\{G(X(n))\}$ (see p.223 of \citet{beran:2013:long}).
\begin{Lem}\label{Lem:cov Herm rank}
If $G(\cdot)$ has Hermite rank $k$, then
\begin{equation}\label{eq:Cov G}
\Cov[G(X(n)),G(X(0))]\approx \Cov[X(n),X(0)]^k \approx n^{(2H-2)k} .
\end{equation}
\end{Lem}
\begin{Rem}\label{Rem:same hurst index}
Comparing (\ref{eq:Cov G}) and (\ref{eq:cov decay}), we note that for functions $G(\cdot)$ that have Hermite rank $k=1$, the Hurst index  of $\{G(X(n))\}$ is the same as the Hurst index  of $\{X(n)\}$. In general  since $2H-2<0$, the higher the Hermite rank $k$ is, the faster the covariance decays as $n\rightarrow\infty$.
Note that in view of (\ref{eq:cov decay}) and (\ref{eq:Cov G}), for $\{G(X(n))\}$ to have long memory, one needs
\[(2H-2)k>-1\iff H>1-{1\over 2k}.
\]
This is natural because when $k>1$, the covariance of $\{G(X(n))\}$ decays faster than that of $\{X(n)\}$ and thus $H$ must be greater than $1-1/(2k)$ in order to ensure that $\{G(X_n)\}$ has long memory.
\end{Rem}

 \smallskip
 \subsection{Asymptotic behavior}
 Now returning to the theme of the introduction: suppose that in order to derive the asymptotic distribution of the statistics $T$ of interest, one first needs to obtain the distributional limit   as $N\rightarrow\infty$ of
\begin{equation}\label{eq:sum stat}
\frac{1}{A(N)}\sum_{n=1}^{[Nt]} \Big( G(X(n))-\E G(X(n)) \Big), \quad t\in [0,1],
\end{equation}
where $G(\cdot)\in L^2(\gamma)$,  $A(N)$ is a suitable normalization, and $[\cdot]$ stands for the integer part.

\begin{Thm}\label{Thm:nclt gaussian}(\citet{dobrushin:major:1979:non},  \citet{taqqu:1979:convergence}, \citet{breuer:major:1983:central}, \citet{major:2014:multiple}.)
Suppose that $G$ has Hermite rank $k$.  Then the following conclusions hold.

\medskip
\noindent $\bullet$   Central limit case: suppose that 
$$
H<1-\frac{1}{2k}.
$$ 
 Then $\{G(X(n))\}$ has short memory in the sense that   \[\sigma^2:=\sum_{n=-\infty}^\infty \Cov[G(X(n)),G(X(0))]\] converges absolutely and
\begin{equation}\label{eq:Herm limit}
\frac{1}{N^{1/2}}\sum_{n=1}^{[Nt]} \Big( G(X(n))-\E G(X(n)) \Big) \ConvFDD \sigma B(t),\quad t\ge 0,
\end{equation}
where $\ConvFDD$ denotes convergence of the finite-dimensional distributions and $B(t)$ is the standard Brownian motion.

\smallskip
\noindent $\bullet$ Non-central limit case:  suppose that 
$$
H>1-\frac{1}{2k}.
$$
Then
$\{G(X(n))\}$ has long memory with Hurst index:
\begin{equation}\label{eq:H_G general}
H_G=(H-1)k+1\in \left(\frac{1}{2},1\right).
\end{equation}
Furthermore, as $N\rightarrow\infty$, we have\footnote{In fact,  we have weak convergence in the space $D[0,1]$  with uniform metric.}
\begin{equation}\label{eq:Herm limit}
\frac{1}{N^{H_G}}\sum_{n=1}^{[Nt]} \Big( G(X(n))-\E G(X(n)) \Big) \ConvFDD  c Z_{H_G,k}(t),
\end{equation}
for some $c\neq 0$,    and
\begin{equation}\label{eq:Herm process}
Z_{H_G,k}(t)=\int_{\mathbb{R}^k}' \left[ \int_0^t \prod_{j=1}^k (s-x_j)_+^{\gamma} ds\right]  B(dx_1)\ldots B(dx_k),\quad \gamma=H-\frac{3}{2}=\frac{H_G-1}{k}-\frac{1}{2},
\end{equation}
is the so-called $k$-th order Hermite process,
where
 $\int_{\mathbb{R}^k}' \left[~  \cdot ~ \right]  \ B(dx_1)\ldots B(dx_k)$ denotes the $k$-tuple Wiener-It\^o integral with respect to the standard Brownian motion $B(\cdot)$.  The prime $'$ indicates that one does not integrate on the diagonals $x_i=x_j$.
\end{Thm}

\begin{Rem}\label{eq:herm rank 1}
When the Hermite rank $k=1$, one has
\[
H_G=H
\]
and the limit $Z_{H_G,k}(t)$ in Theorem \ref{Thm:nclt gaussian}  is  the \emph{fractional Brownian motion} $B_{H}(t)$, namely, the centered Gaussian process determined by the following covariance structure:
\[
\Cov\left[B_{H}(s),B_{H}(t)\right]=\frac{1}{2}\left( |s|^{2H}+|t|^{2H}-|s-t|^{2H} \right).
\]
The preceding covariance is shared by all the other Hermite processes.
When the Hermite rank $k=2$, $Z_{H_G,2}(t)$ is called the \emph{Rosenblatt process} (see \citet{rosenblatt:1961:independence} and \citet{taqqu:1975:weak}).
The Hermite process  $Z_{H_G,k}(t)$ in (\ref{eq:Herm process}) admits different representations. See \citet{pipiras:taqqu:2010:regularization}.
\end{Rem}
\begin{Rem}\label{Rem:breuer major}
  The boundary case $H=1-\frac{1}{2k}$ typically  falls in the central limit theorem  regime (convergence to Brownian motion) after modifying the normalization $N^{-1/2}$ to include some slowly varying functions  (Theorem 1' of \citet{breuer:major:1983:central}).  In general, the convergence of finite-dimensional distributions $\ConvFDD$ in the short-memory case   cannot be strengthened to weak convergence $\Rightarrow$ in $D[0,1]$ unless some additional assumption is imposed on $G$, e.g.,  $G$ being a polynomial (\citet{chambers:slud:1989:central}).
\end{Rem}

The long-memory Gaussian  $\{X(n)\}$   may be directly used as a model for the long-memory stationary noise. For statistical theory, however, it is  often desirable to allow  departure from   Gaussianity, e.g., to accommodate the situations where the noise distribution is  skewed or heavy-tailed.   Within the same framework,  a way to achieve such  flexibility is as follows.
Suppose  that there is an underlying long-memory \emph{Gaussian}  stationary process $\{Y(n)\}$. Assume without loss of generality  that $\{Y(n)\}$ is standardized.  Now suppose  that the noise sequence $\{X(n)\}$ in the model  is given by
\begin{equation}\label{eq:inst trans}
X(n)=F(Y(n)).
\end{equation}

\begin{Rem}
There are different perspectives to interpret  (\ref{eq:inst trans}). First, note that
when $F(\cdot)$ is nonlinear, $X(n)$ is non-Gaussian. So $F(\cdot)$ can represent the  departure from the ideal Gaussian assumption. Hence when the noise $X(n)$ is  modeled by (\ref{eq:inst trans}) with  an unknown $F(\cdot)$, this provides great model flexibility. Note that a proper choice of $F(\cdot)$ can match any marginal distribution for $X(n)$. Second, from the perspective of analysis of robustness, one may  view $X(n)$ as a perturbed version of $Y(n)$, where $X(n)$ is close to $Y(n)$, that is, $F(\cdot)$ is close to the identity function.
\end{Rem}

Following the same statistical inference procedure that leads to (\ref{eq:sum stat}),  we then focus on  the distributional limit  of
\[
 \frac{1}{A(N)}\sum_{n=1}^{[Nt]} \Big( G\circ F (Y(n))-\E G\circ F(Y(n)) \Big)
 \]
as $N\rightarrow\infty$.
\begin{Rem}
We emphasize the different roles played by $F(\cdot)$ and $G(\cdot)$. \textbf{The function $F(\cdot)$   accounts for  an unknown and uncontrollable departure from the Gaussian $Y(n)$.    On the other hand, the function $G(\cdot)$ depends on the statistical procedure of interest and is therefore typically precisely known}. For example,   $G$  is typically the identity transformation for inference of the mean $\E X(n)=\E F(Y(n))$.
\end{Rem}

\smallskip
\subsection{Basic claim}
We are now ready   to make the following  claim which will be justified below. The case of  non-instantaneous (multivariate) $F$ (and also $G$)   will be addressed in Appendix \ref{Sec:non-inst trans}   (the issues remain essentially the same).
\begin{Cla}\label{Cla:Herm}
It is typically the case that
\begin{enumerate}[(a)]
\item the function $G\circ F$ has Hermite rank 1;
\item the process $\{X(n)\}$ has long memory with the same Hurst index $H$ as $\{Y(n)\}$ in (\ref{eq:inst trans}).
\end{enumerate}
\end{Cla}

\medskip
\begin{proof}[Justification of the claim:]

Let  $Z$ be  standard Gaussian. Since $H_1(Z)=Z$, requiring   the function $G\circ F$ to have Hermite rank $k\ge 2$ is equivalent to
\begin{equation}\label{eq:herm rank >= 2}
\E[ (G\circ F)(Z)Z]=0.
\end{equation}
This requirement is very restrictive, and   is, moreover, unrelated to the usual size or smoothness conditions typically imposed on the perturbation $F$.
\textbf{\emph{Unlike the precisely known $\boldsymbol{G(\cdot)}$ which is related to   the method of   inference considered, one has no control  nor accurate knowledge of the  function $\boldsymbol{F(\cdot)}$}}.
 There is thus no a priori reason     that $F(\cdot)$ be  such that (\ref{eq:herm rank >= 2}) holds. But if (\ref{eq:herm rank >= 2}) does not hold, then the Hermite rank of $G\circ F$ is $1$, which justifies part (a) of the Claim \ref{Cla:Herm}.

Applying the same  reasoning, the perturbation function $F(\cdot)$ is also very likely to be such that
\[\E[ F(Z)Z]\neq 0
\] and hence  to have Hermite rank $1$. Then in view of  Lemma \ref{Lem:cov Herm rank} and Remark \ref{Rem:same hurst index}, this justifies part (b) of the Claim \ref{Cla:Herm} (b).
\end{proof}
\begin{Rem}
The Claim \ref{Cla:Herm} indicates not only  the instability of a Hermite rank higher than $1$, but also the stability of the Hermite rank 1 and hence the  Hurst index of the noise model.  Then, as suggested by the Claim \ref{Cla:Herm},  if $G\circ F$ has Hermite rank 1,  by Theorem \ref{Thm:nclt gaussian} and Remark \ref{eq:herm rank 1}, one has
\begin{equation}\label{eq:conv to fBm}
\frac{1}{N^H}\sum_{n=1}^{[Nt]} \Big( G\circ F (Y(n))-\E  G\circ F (Y(n)) \Big) \Rightarrow cZ_{H,1}(t)  = c B_H(t).
\end{equation}
for some $c>0$, where $B_H(t)$ is the fractional Brownian motion with Hurst index $H$. The theorem thus also implies the stability of  fractional Brownian motion as the limit.
\end{Rem}

\begin{Rem}
In statistics, one sometimes needs limit theorems for  functionals other than the sum. A typical example is the quadratic form $\sum_{n,m=1}^{N} a(n-m)  X(n) X(m)$. Limit theorems in this case depends on the not only the ``memory'' of $\{X(n)\}$ but also the ``memory'' of the coefficient $a(n)$ (see, e.g., \citet{avram:1988:bilinear} and \citet{terrin:taqqu:1990:noncentral}). Instead of discussing in general the instability of such     quadratic forms,   we shall focus in Section \ref{Sec:whittle} below on an important statistical application, namely,  Whittle estimation.
\end{Rem}


\smallskip
\subsection{The level shift case}

One may consider making Claim \ref{Cla:Herm} a genuine mathematical statement by, for example, considering $F(\cdot)$ as a \emph{random} element in a suitable function space with a ``prior probability model'', as long as that   model assigns a small probability to the set of $F(\cdot)$ on which (\ref{eq:herm rank >= 2}) happens.

  In the following theorem,
we consider the simple case  where the perturbation is given by  a level shift of size $z$, namely, if $F(y)=z+y$  so that $G\circ F(\cdot)=G(\cdot+z)$. To understand the assumptions, note that we want to exclude the case where $G(\cdot)$ is constant, since then $G(\cdot+z)$ remains equal to $G(\cdot)$. We also want $G(\cdot+z)$ to be in $L^2(\gamma)$.
\begin{Thm}\label{Thm:shift}
Suppose that the function $G(\cdot)\in L^2(\gamma)$   has an arbitrary Hermite rank,
  $G(\cdot)$ is not  constant a.e.\ , and
assume that there exists $\delta>0$, so that $G(\cdot + z)\in L^2(\gamma)$ for all $|z|<\delta$. Then there exists $\epsilon\in (0,\delta)$, such that the Hermite rank of $G(\cdot +z )$ is $1$ for all $z\in (-\epsilon,0)\cup (0,+\epsilon)$.
\end{Thm}
 The proof can be found in \citet{bai:taqqu:2017:hermite}. In that paper, we also study what happens when    the shift tends to zero as the sample size tends to infinity, which is analogous to  the near integration analysis of unit roots (see \citet{phillips:1987:towards}). In   \citet{bai:taqqu:2017:hermite}, we also consider    transformations other than the shift, e.g., the scaling   $F(z)=zy$ so that $G\circ F(y)=G(zy)$.

\smallskip
\subsection{Transformation of linear processes}\label{Sec:linear}
Another popular class of  models for a stationary, not necessarily Gaussian,  noise $\{Y(n)\}$  is the so-called (causal) linear process:
\begin{equation}\label{eq:linear}
Y(n)=\sum_{i=0}^\infty a_{n-i} \epsilon_i,
\end{equation}
where $\epsilon_i$'s are assumed to be i.i.d.\  random variables (not necessarily Gaussian) with mean $0$ and  variance $1$ and $\sum_n a_n^2<\infty$.
When
\begin{equation}\label{eq:a_n long memory}
a_n\approx n^{H-3/2},\quad 1/2<H<1 \quad\text{ as } \quad n\rightarrow\infty,
\end{equation}
one has $\Cov[Y(n),Y(0)]\approx n^{2H-2}$, and thus $Y(n)$ has long memory with Hurst index $H$. The well-known fractionally-integrated noise model (see, e.g., \citet{granger:joyeux:1980introduction}) satisfies  (\ref{eq:a_n long memory}).
 We shall assume (\ref{eq:a_n long memory}) throughout this section.

Theorem \ref{Thm:nclt gaussian}  can be  extended to  linear processes. In this case,   the larger class of polynomials called   Appell polynomials (\citet{avram:1987:noncentral}) plays an analogous role to that of the Hermite polynomial in Section \ref{Sec:inst trans}. One can define the so-called \emph{Appell} rank  of a function $G(\cdot)$ as in (\ref{eq:def herm rank expan}),  with Hermite polynomials  replaced by the Appell polynomials, given that the expansion is valid (for example, when $G(\cdot)$ is a finite-order polynomial). However, in this    framework (\citet{surgailis:1982:zones}), the class of functions $G(\cdot)$'s that can be treated is rather restrictive. \citet{ho:hsing:1997:limit}  greatly extended   the allowable $G(\cdot)$'s through a martingale difference approach and  introduced a more convenient notion of rank, which we shall call the \emph{power rank}.  See also \citet{levy-leduc:taqqu:2013:long}.

Given a function $G(\cdot)$ and a random variable $Y$ satisfying $\E G(Y)^2<\infty$, let
\begin{equation}\label{eq:G_infty}
G_\infty(y)=\E  G (Y+y)
\end{equation}
given that the expectation exists  and suppose that $G_\infty(\cdot)$ has   derivatives of  order sufficiently high. The power rank of $G(\cdot)$ with respect to $Y$ is defined as
\begin{equation}\label{eq:power rank}
\inf\{m \ge 1:~ G_\infty^{(m)}(0)\neq 0  \},
\end{equation}
where $G_\infty^{(m)}(y)$ denotes the $m$-th derivative of $G_\infty(y)$. In fact, the power rank in (\ref{eq:power rank}) coincides with the Hermite rank \eqref{eq:def herm rank} if $Y$ is Gaussian.
This was stated   in \citet{ho:hsing:1997:limit}, and see \citet{bai:taqqu:2017:hermite} for a proof.
The case  where $G$ is a polynomial was treated in \citet{levy-leduc:taqqu:2013:long}.

Now we can state the following  limit theorems (for simplicity we omit  the inclusion of some technical conditions, see \citet{ho:hsing:1997:limit} and \citet{pipiras:taqqu:2017:long}):
\begin{Thm}[\citet{ho:hsing:1997:limit}]\label{Thm:nclt linear}
Suppose that $\E G(Y(0))^2<\infty$ and $G(\cdot)$ has power rank $k\ge 1$ with respect to $Y(0)$ in the sense of (\ref{eq:power rank}).
Under some additional technical conditions, statements exactly analogous  to Theorem \ref{Thm:nclt gaussian} hold  with the role of the Hermite rank replaced by the power rank.
\end{Thm}

\begin{Rem}
Using similar arguments as Section \ref{Sec:inst trans}, one sees that a power rank higher than $1$ is also unstable to  perturbation: to get a power rank higher than $1$, one needs the restrictive condition
\[G_\infty^{(1)}(0)=\frac{d}{dy}\E G(Y+y)|_{y=0}=0.
\]
which can be easily perturbed by compositing $G$ with a transformation before. So an   analog of Claim \ref{Cla:Herm} may   be stated in this context.
\end{Rem}

Below we provide some further  remarks on the instability phenomenon.


\begin{Rem}
We mention that
\citet{surgailis:2000:long} established some results which can be interpreted as the ``robustness'' of Theorem \ref{Thm:nclt linear} against additive noise.  Roughly speaking, \citet{surgailis:2000:long} showed that  if the long memory linear process $Y(n)$ is replaced by $Y'(n)=Y(n)+U(n)$ with $U(n)$ specified as some short-memory models, then the non-central convergence in Theorem \ref{Thm:nclt linear} still holds, where the power (or Appell) rank  is now with respect to the distribution of $Y'(n)$. Nevertheless, the instability discussed earlier still applies.
   Firstly, the rank can still be unstable under a  transformation.
Secondly,  even without considering a   transformation perturbation,   one has typically no accurate knowledge of the marginal distribution of $Y'(n)$.  A change in   the distribution  will affect the  rank  defined through (\ref{eq:power rank}).
\end{Rem}


\begin{Rem}
It is important to consider  not only the limit distribution that one obtains, but also  the normalization since the latter  corresponds to the magnitude of the fluctuations of the partial sum. When the true rank (Hermite or power) is indeed $1$, but one assumes a higher-order Hermite rank from some statistical consideration,  this  will lead to underestimation of the magnitude of the fluctuations of the partial sum since $H_G<H$ in Theorem \ref{Thm:nclt gaussian}.
\end{Rem}

\section{Examples in Statistics}\label{Sec:examples}

In this section, we  review  statistical problems in   the  literature  related to  limit theorems involving different ranks.  We shall elaborate on some examples: sample variance, empirical processes, Whittle likelihood, and nonparametric estimation, to demonstrate how the asymptotic statistical theories  are affected  by the instability  discussed in Section \ref{Sec:Instability}.

\smallskip
\subsection{Sample variance}\label{Sec:var}

In the context of long memory, with the complexity introduced by  the limit theorems,   scale estimation becomes  a problem. We only discuss the case where the data $\{X(n)\}$ is a Gaussian process, but everything   can be extended to a linear process  $\{X(n)\}$  (see Section \ref{Sec:linear}).

Assume then that $\{X(n)\}$ is a long memory stationary Gaussian process with Hurst index $H\in (1/2,1)$, unknown mean $\mu$ and unknown variance $\sigma^2$. Consider the estimation of $\sigma
^2$ using the sample variance:
\begin{equation}
\widehat{\sigma}_N^2:=\frac{1}{N}\sum_{n=1}^N \left(X(n)-\bar{X}_N\right)^2,\label{e:GX}
\end{equation}
where $\bar{X}_N=(X(1)+\ldots+X(N))/N$ is the sample mean.  In the short memory situation, say if $X(n)$ were i.i.d., it is well known that $\widehat{\sigma}_N^2$ is asymptotically normal. The situation is, however, delicate.
Indeed,  express (\ref{e:GX}) as
\begin{equation}\label{eq:sample var}
\widehat{\sigma}_N^2=\frac{1}{N}\sum_{n=1}^N (X(n)-\mu)^2+  (\bar{X}_N-\mu)^2=: U_N+V_N.
\end{equation}
We can write
\begin{equation}\label{eq:U_N}
U_N-\sigma^2=\frac{1}{N}\sum_{n=1}^N (X(n)-\mu)^2-\sigma^2= N^{-1/2} \left[ \frac{1}{\sqrt{N}} \sum_{n=1}^N \left((X(n)-\mu)^2-\sigma^2\right)\right],
\end{equation}
and
\begin{equation}\label{eq:V_N}
V_N=  (\bar{X}_N-\mu)^2=N^{2H-2}\left[ \left( \frac{1}{N^{H}}\sum_{n=1}^N (X(n)-\mu)\right)^2\right].
\end{equation}
Note that the term $X(n)-\mu$ has Hermite rank $k=1$,
and the term $(X(n)-\mu)^2$ has    expectation $\sigma^2$ and Hermite rank $k=2$, since $\E (X(n)-\mu)^2 X(n)=0$.

Thus when $H<3/4$, in view of Theorem \ref{Thm:nclt gaussian},   the term in the brackets in the right-hand sides (\ref{eq:U_N}) and (\ref{eq:V_N}) converge as $N\rightarrow\infty$.
Since $H<3/4$ implies $N^{2H-2}\ll N^{-1/2}$, the term $V_N$ is asymptotically negligible, compared to $U_N-\sigma^2$, and hence $N^{1/2}( \widehat{\sigma}_N^2-\sigma^2)$ has the same limit as $N^{1/2}U_N$ as $N\rightarrow\infty$. Thus
\[
N^{1/2}( \widehat{\sigma}_N^2-\sigma^2) \ConvD N(0,s_1^2)
\]
for some $s_1>0$.

 When $H>3/4$,  in view of Theorem \ref{Thm:nclt gaussian} with $k=2$, we write
\[
U_N-\sigma^2= N^{2H-2} \left[\frac{1}{N^{2H-1}} \sum_{n=1}^N \left((X(n)-\mu)^2-\sigma^2\right)\right],
\]
and $V_N$ is as in (\ref{eq:V_N}). Now both $U_N$ and $V_N$ contribute to the limit, where we have by a multivariate version of Theorem \ref{Thm:nclt gaussian} (see, e.g, \citet{bai:taqqu:2013:1-multivariate}):
\begin{equation}\label{eq:limit mix}
N^{2-2H}( \widehat{\sigma}_N^2-\sigma^2) \ConvD a_H Z_{2H-1,2}(1)+  b_H Z_{H,1}(1)^2
\end{equation}
where $Z_{2H-1,k}(t), k=1,2$ are the Hermite processes in (\ref{eq:Herm process}) defined by the same Brownian integrator $B(\cdot)$.
for some constants $a_H,b_H>0$. See also \citet{dehling:taqqu:1991:bivariate}.

We now suppose that $\{X(n)\}$ is perturbed by a  transformation in the spirit of Claim \ref{Cla:Herm}, which leads to consider the case where both $X(n)$ and $[X(n)-\mu]^2$  have Hermite rank $1$.  Then  writing
\[
U_N-\sigma^2= N^{H-1} \left[\frac{1}{N^{H}} \sum_{n=1}^N \left((X(n)-\mu)^2-\sigma^2\right)\right],
\]
and
\begin{equation*}
V_N= N^{2H-2} \left[ \frac{1}{N^{H}}\sum_{n=1}^N (X(n)-\mu)\right]^2,
\end{equation*}
we can apply Theorem \ref{Thm:nclt gaussian} with $k=1$.
Since $H<1$, we have $N^{2H-2}\ll N^{H-1}$, and thus only the term $U_N$ contributes to the limit.
Then $N^{1-H}(\widehat{\sigma}_N^2-\sigma^2)$ has the same limit as
\[N^{1-H}(U_N-\sigma^2)=\frac{1}{N^{H}}\sum_{n=1}^N [(X(n)-\mu)^2-\sigma^2],
\]
namely, $c_H Z_{H,1}(1)$
for some $c_H>0$, where $Z_{H,1}$ is the fractional Brownian motion. Hence
\begin{equation}\label{e:Nsig}
N^{1-H} \left(\widehat{\sigma}_N^2-\sigma^2\right) \ConvD c_H Z_{H,1}(1),
\end{equation}
which is different from (\ref{eq:limit mix}).

\begin{Rem}
Under the above perturbation consideration, there is no dichotomy between $H<3/4$ and $H>3/4$ in (\ref{e:Nsig}), and the normalization  is always $N^{1-H}$, which is of smaller order than  both $N^{1/2}$ and $N^{2H-2}$. In the case $H<3/4$, however,  we get a Gaussian limit with or without perturbation. Hence   without the perturbation consideration, there is  the danger of underestimating the fluctuation magnitude of the sample variance, namely, taking the fluctuation to be of the order $N^{-1/2}$ when $H<3/4$ and $N^{2H-2}$ when $H>3/4$, whereas they are of the order $N^{H-1}$.  We also mention that  similar considerations also  apply to the study of the asymptotic behavior of sample autocovariance/correlation (see, e.g., \citet{hosking:1996:asymptotic},  \citet{wu:huang:2010:covariance} and \citet{levy-leduc:2010:robust}).
\end{Rem}

\smallskip
\subsection{Empirical processes}
Empirical processes play important roles in many statistical problems. We refer the reader to    \citet{dehling:mikosch:2002:empirical} for an introduction to empirical processes of dependent data. Let $\{X(n)\}$ be a stationary process with marginal cdf $F(x)$. The corresponding centered empirical process is defined as
\begin{align}\label{eq:empirical}
F_N(x)=\frac{1}{N}\sum_{n=1}^{N}\left[ I\{X(n)\le x\}-F(x) \right].
\end{align}
When $\{X(n)\}$ is i.i.d., it is well-known that $N^{1/2}F_N(x)$ converges weakly in $D(-\infty,\infty)$ to  $B_0(F(x))=B(F(x))-F(x)B(1)$, where $B_0(t)=B(t)-tB(1)$  is a Brownian bridge and $B(t)$ is a Brownian motion. Under some weak dependence conditions  on $\{X(n)\}$, the process $N^{1/2}F_N(x)$ converges weakly  in $D(-\infty,\infty)$ to a centered  Gaussian process $G(x)$ with covariance structure given by
\begin{equation}\label{eq:emp cov}
\E G(x)G(y)=\sum_{n=-\infty}^\infty \Cov[ I(X(0)\le x),I(X(n)\le y) ].
\end{equation}
See, e.g., Theorem 4.1 of \citet{dehling:philipp:2002:empirical}.

When $X(n)$ has long memory, the corresponding weak convergence results become rather different in nature. Indeed, assume that $X(n)=G(Y(n))$ where $\{Y(n)\}$ is a standardized stationary Gaussian process with Hurst index $1/2<H<1$. We define the deterministic function
\[
J_m(x)=\frac{1}{m!}\E I\{G(Y(0))\le x\}H_m(Y(0)),
\]
where $H_m(\cdot)$ is the $m$-th order Hermite polynomial. Note that for any fixed $x\in \mathbb{R}$, $J_m(x)$'s are the coefficients of the Hermite expansion (\ref{eq:def herm rank expan}) of the function $\Delta_x(y)=I\{G(y)\le x\}-F(x).$
We have the following result:
\begin{Thm}[Theorem 1.1 of \citet{dehling:taqqu:1989:empirical}]\label{Thm:empirical}
Let
\begin{equation}\label{eq:empirical rank}
k=\inf\{m\ge 1:~J_m(x)\neq 0 \text{ for at least one }x\in\mathbb{R}\},
\end{equation}
and assume that $H>1-\frac{1}{2k}$. Then we have the following weak convergence in $D(-\infty,+\infty)$
\begin{equation}\label{eq:empirical conv}
N^{1-H_k} F_N(\cdot) \Rightarrow c J_m(\cdot)Z_{H_k,k}(1)
\end{equation}
where $Z_{H_k,k}(\cdot)$ is the Hermite process as in (\ref{eq:Herm process}), and $H_k=(H-1)k+1$  as in (\ref{eq:H_G general}).
\end{Thm}
\begin{Rem}
It is interesting to note that in the long memory case, the limit process  
$$ \{J_m(x)Z_{H_k,k}(1),~x\in \mathbb{R}\}$$
is quite degenerate, namely, it has correlation $1$ between any different points $x_1,x_2\in \mathbb{R}$, in contrast to the weak dependence  case where the limit Gaussian process $G(x)$ admits a rich correlation structure (see (\ref{eq:emp cov})).
\end{Rem}

By the perturbation argument,   one may assume that the rank $k=1$, regardless of the choice of $G(\cdot)$ in a statistical application of Theorem \ref{Thm:empirical}. In fact, the definition of rank  (\ref{eq:empirical rank})  makes the assumption $k=1$ even more appealing in this context, because $J_1(x)\neq 0$ for just one point $x$ would make $k=1$. From this point of view, the only practically relevant convergence in (\ref{eq:empirical conv}) is
\begin{equation}\label{eq:emp conv k=1}
N^{1-H} F_N(\cdot) \Rightarrow cZ  J_1(\cdot),
\end{equation}
where $Z$ is a standard Gaussian variable, and thus the fluctuation of the empirical process is practically always of the order $N^{H-1}$. The convergence (\ref{eq:empirical conv}) can be applied to study the asymptotic behavior of U-statistics and V-statistics (see Corollary 1 of \citet{dehling:taqqu:1989:empirical}).   It is also applied to develop the asymptotic theories of estimation of the  probability density function $f=F'$ (see \citet{csorgo:mielniczuk:1995:density} and Section \ref{Sec:misc}).

\smallskip
\subsection{Whittle likelihood}\label{Sec:whittle}
In the parametric estimation for time series, the so-called Whittle pseudo-likelihood  is a  computationally efficient approximation to the Gaussian likelihood, which bypasses the inversion of a covariance matrix in the latter.
The resulting Whittle estimator and its semiparametric extensions are found particularly useful  in the long memory context for  the estimation of the Hurst parameter $H$.  For more details on the background and motivation, we refer  to Section 5.5 of \citet{beran:2013:long} or Chapter 10 of \citet{pipiras:taqqu:2017:long}. We shall focus on the rank instability issue in the asymptotic theory developed  in \citet{giraitis:taqqu:1999:whittle}. The asymptotic theory in this context depends on the limit theorem for quadratic forms,  which is more delicate than the limit theorems for sums. The instability  issue in this context exhibits some distinct features compared with the previous cases.

Suppose a stationary time series $\{X(n)\}$ has spectral density (see, e.g., Chapter 1 of \citet{pipiras:taqqu:2017:long})   $f(\lambda;\theta,\sigma)=\sigma^2g_\theta(\lambda)>0$, $\lambda\in (-\pi,\pi)$ so that
\[
\Cov[X(n),X(0)]=\int_{-\pi}^{\pi} e^{inx} f_\theta(\lambda;\theta,\sigma)d\lambda,
\]
where $\sigma$ and $\theta=(\theta_1,\ldots,\theta_p)$ are unknown parameters.
Assume that the normalization  condition (scale identifiability) holds:
\[
\int_{-\pi}^{\pi} \log g_\theta(\lambda) d\lambda=0,
\]
 under which $\sigma^2$ becomes the mean squared error of the one-step prediction  by the Kolmogorov's formula (see, e.g., Section 5.8 of \citet{brockwell:1991:time}).
Suppose that we want to estimate the unknown parameter $\theta$. Under long memory, the choice of $\theta$   typically includes $H$.  Define
\begin{equation}\label{eq:a_theta}
a_\theta(n)= \int_{-\pi}^{\pi} e^{in\lambda} \frac{1}{g_\theta(\lambda)} d\lambda
\end{equation}
and
\[
w_\theta(m,n)=\int_{-\pi}^{\pi} g_\theta(\lambda) \frac{\partial^{2}}{\partial \theta_m \partial \theta_n} [g_\theta(\lambda)]^{-1} d\lambda.
\]
The so-called Whittle estimator $\widehat{\theta}_N$ of $\theta$ is given by
\[
\widehat{\theta}_N=\underset{\theta}{\mathrm{argmin}} \sum_{m,n=1}^N  a_\theta(m-n) X(m) X(n).
\]
If $\{X(n)\}$ is a Gaussian or a linear long-memory process,  it was  established under some regularity conditions that (see, e.g., \citet{fox:taqqu:1986:large} and \citet{giraitis:surgailis:1990:central})
\begin{equation}\label{eq:Whittle CLT}
N^{1/2} (\widehat{\theta}_N-\theta)\ConvD  N(0,4\pi W^{-1}_\theta),
\end{equation}
where the matrix $W_\theta=(w_\theta(m,n))_{1\le m,n\le p}$. Note that the standard $N^{1/2}$-convergence rate appears even though $\{X(n)\}$ has long memory.  This is  due to the dependence cancellation  effect from the quadratic coefficient $ a_\theta(n)$.

On the other hand,  \citet{giraitis:taqqu:1999:whittle}, considered
$
X(n)=G(Y(n))$,
where $\{Y(n)\}$ is long-memory Gaussian  and the transformation $G(\cdot)$ is restricted to be a polynomial by \citet{giraitis:taqqu:1999:whittle} to avoid some technical difficulties.
Define
\begin{align*}
\rho_k&=\frac{1}{k!} \sum_{n=-\infty}^{\infty} \E\left[ \frac{d^k}{dx^k}G(x+Y(0))G(x+Y(n))\right]\biggr\rvert_{x=0} \nabla a_{\theta}(n)
\\& =\sum_{m,n\ge 0, m+n=k} \frac{1}{m!n!} \sum_{r=-\infty}^\infty \E \left[G^{(m)}(Y(r))  G^{(n)}(Y(0)) \right] \nabla a_{\theta}(r),
\end{align*}
where  $\nabla$ denotes the gradient with respect to $\theta$. In particular,
\begin{equation}\label{eq:rho_1}
\rho_1=2 \sum_{n=-\infty}^\infty \E[ G' (Y(n)) G(Y(n))]\nabla a_\theta(n).
\end{equation}
Note that in the case $G(x)=x$, namely, the Gaussian case, $\rho_k=0$ for all $k=1,2,\ldots$.

It was established in Corollary 2.1 of \citet{giraitis:taqqu:1999:whittle} that  under some regularity conditions, if $\rho_1\neq 0$, then  as $N\rightarrow\infty$, we have
\begin{equation}\label{eq:Whittle rho_1}
N^{1-H} (\widehat{\theta}_N-\theta) \ConvD  Z
\end{equation}
for some centered normal random vector $Z$. Note that  in (\ref{eq:Whittle rho_1}) the convergence rate is the same as that of the sample mean in view of (\ref{eq:var allan}).
\citet{giraitis:taqqu:1999:whittle} also showed that if $\rho_1=0$ but some $\rho_k\neq 0$, then (\ref{eq:Whittle rho_1}) needs to be modified resulting in limit theorem with a central and non-central dichotomy similar to Theorem \ref{Thm:nclt gaussian}. See Theorem 2.3 and 3.1 of \citet{giraitis:taqqu:1999:whittle}.

Now we consider the instability issue. Here the role of Hermite (or power) rank is instead played by
\begin{equation}\label{eq:rank whittle}
k=\inf\{m\in \mathbb{Z}_+:~\rho_m\neq 0\}.
\end{equation}
  There is instability even in the Gaussian case  where $G(\cdot)$ is the identity, namely, $G(x)=x$, and where then all $\rho_k= 0$. In that case, we would have (\ref{eq:Whittle CLT}), but by perturbing $G(\cdot)$ slightly, we would get $\rho_1\neq 0$ in (\ref{eq:rho_1}), and thus we would have (\ref{eq:Whittle rho_1}) instead of (\ref{eq:Whittle CLT}).

The preceding observation raises a question on  the applicability of (\ref{eq:Whittle CLT}) in statistical inference. It turns out that the achievement of the parametric rate $N^{1/2}$, or say  the cancellation effect  of the quadratic coefficient in (\ref{eq:a_theta}),  critically depends on the Gaussian or linear data-generating assumption, while a disturbance of such an assumption yields instead the rate $N^{1-H}$, which is the usual slower rate  of convergence under long memory.
It is unclear whether similar instability issues occur in  the semiparametric extensions of the Whittle estimator, e.g., the local Whittle estimator (\citet{kunsch:1987:statistical}, \citet{robinson:1995:gaussian}).

\smallskip
\subsection{Nonparametric estimation}\label{Sec:nonpara}
In this section, we  review briefly some nonparametric statistical studies under long memory  involving  the  Hermite or power rank.    Assume  throughout that $\{X(n)\}$ is a stationary \emph{long-memory} process, typically specified by a Gaussian process, or a linear process, or a  transformation of either (we call the model a Gaussian or linear subordination).

In the kernel smoother type nonparametric estimation procedures, a nonlinear transformation of the data is naturally involved.
For example, the kernel density estimator of the probability density function $f(x)$ is defined as
\begin{equation}\label{eq:density est}
\widehat{f}(x)= \frac{1}{Nh}\sum_{n=1}^N  K\left( \frac{x-X(n)}{h}\right), \qquad x\in \mathbb{R}.
\end{equation}
where $N$ is the sample size, $h>0$ is the bandwidth parameter, and $K(\cdot)$ is a kernel satisfying $\int_{\mathbb{R}}K(x)dx=1$. A number of studies have considered the asymptotic behavior of the estimator $\widehat{f}(x)$ as $N\rightarrow\infty$ and $h\rightarrow 0$. See, for example, \citet{cheng:robinson:1991:density}, \citet{csorgo:mielniczuk:1995:density} and \citet{ho:1996:central},\citet{wu:mielniczuk:2002:kernel}.
Another typical class of statistical procedures involving kernel smoothers are the nonparametric local regressions  (e.g., Nadaraya-Watson estimator and local polynomial estimators). Some relevant work involving  the  ranks are \citet{hidalgo:1997:non} , \citet{csorgo:mielniczuk:1999:random},  \citet{masry:mielniczuk:1999:local}, \citet{guo:koul:2007:nonparametric}.

We  mention that under long memory, the asymptotic behavior of kernel smoothers can be quite different from the short-memory case. In particular, an interesting dichotomy phenomenon appears in the asymptotics depending on how fast the bandwidth $b_n$ tends to $0$ with respect to the Hurst index $H$ of $X(n)$. If $b_n$ tends to $0$ relatively slowly, one can have a very degenerate behavior such as the kernel density estimate   $\widehat{f}(x)$ at different points of $x$ becomes  asymptotically perfectly correlated.  See, e.g, \citet{csorgo:2002:smoothing} as well as Chapter 5.14 of \citet{beran:2013:long}).

Asymptotic results involving applying limit theorems with different ranks when studying these kernel smoother procedures are due more often to  the   assumption that $X(n)$ is a  transformation of a Gaussian or linear process,  than due to the nonlinear transformation produced by   $K(\cdot)$ in (\ref{eq:density est}).  To obtain a higher-order rank for $\widehat{f}(x)$ when $X(n)$ is Gaussian or linear, one has to be in very special situations, e.g.,  when focusing on the asymptotic distribution of $\widehat{f}(x)$ in (\ref{eq:density est}) at a fixed point $x=x_0$  while assuming that the true density satisfies $f'(x_0)=0$ (see, e.g., Theorem 3 of \citet{wu:mielniczuk:2002:kernel}).
Similar considerations extend to other nonparametric procedures, e.g, the spline regression under long memory  noise (\citet{beran:2011:spline}).

\smallskip
\subsection{Miscellaneous}\label{Sec:misc}

Wavelets are   useful tools for  analyzing long-memory data due to their natural adaptivity to scaling. Limit theorems  involving ranks were applied, e.g.,   in \citet{clausel:2012:large} and \citet{clausel:2014:wavelet}, who studied the asymptotic behaviors of the wavelet coefficients  and the wavelet estimation of Hurst index of the Gaussian subordination data.

Some other statistical studies involving higher-order ranks   limit theorems  are:  bivariate U-processes (\citet{levy:boistard:taqqu:reisen:2011:asymptotic}), change-point test (\citet{zhao:2010:ratio}, \citet{dehling:rooch:2013:non}),  goodness-of-fit test in regression \citet{koul:1998:regression}, normality test (\citet{beran:gosh:1991:slowly}), sign test (\citet{psaradakis:2010:inference}), unit root test  (\citet{wu:2006:unit}).

\section{Empirical evidence}\label{Sec:empirical}
In this section, we provide  empirical evidence to support the preceding discussion of instability of   ranks in the   limit theorems under long memory.

Consider the  rank of the  quadratic transformation 
$$ G(x)=x^2.
$$
 It  is always $2$,  in both the Gaussian and linear subordination context. This means that if  $\{X(n)\}$ is exactly a centered Gaussian or linear process with  Hurst index $H>1/2$, then the Hurst index of the transformed  series $\{X(n)^2\}$ should be
\begin{equation}\label{eq:H_G rank 2}
H_G=\max\left(\frac{1}{2},  2H-1\right)<H,
\end{equation}
 in view of Theorem \ref{Thm:nclt gaussian}   and \ref{Thm:nclt linear}. Note that when $H<0.75$, the resulting Hurst index is always $H_G=0.5$ unless in the special case where the sum of covariances of all orders is zero (ainti-persistency).

Here is the question:  if $\{X(n)\}$ is a  real-life centered stationary data in which displays long memory,    does one typically observe the decrease from $H$ to $H_G$  as in (\ref{eq:H_G rank 2}) when $\{X(n)\}$ is replaced by $\{X(n)^2\}$? If  our  arguments in the previous sections   make practical sense,  then the time series $\{X(n)^2\}$  should most likely still possess rank $1$, which means that (\ref{eq:H_G rank 2}) should barely happen. To test this hypothesis, we design the following empirical study which involves $\{X(n)\}$ and $\{X(n)^2\}$.
The design is explained in Remark \ref{Rem:explain desgin} below.

\bigskip
\noindent{\it
\textbf{Design of the study}:

\bigskip
Suppose that we have a collection of $M$ real-life stationary long-memory time series data
\[\{X_m(n),~n=1,\ldots,N_m~m=1,\ldots,M. \},\]
 where $n$ is the time index, and $m$ is the data set index.
For each $m$, we perform the following analysis.
\begin{enumerate}[\underline{Step }1]
\item For each $m=1,\ldots,M$, center the data:
 $$
 X_m(n)\leftarrow  X_m(n)- \frac{1}{N_m}\sum_{n=1}^{N_m}X_m(n);
 $$
 
 \smallskip
\item  For each $m=1,\ldots,M$, obtain 
\[
\mbox{\it 
the estimated   Hurst index $\widehat{H}_{m}^{(1)}$ of $\{X_m(n)~n=1,\ldots,N_m\}$
}\]
 and 
\[
\mbox{\it  the estimated Hurst index $\widehat{H}_{m}^{(2)}$ of $\{X_m(n)^2,~n=1,\ldots,N_m\}$.
}\]

\smallskip
\item For each $m=1,\ldots,M$, simulate $R$ independent sequences of fractional Gaussian noise (increments of fractional Brownian motion): 
$$
\{G_{mr}(n),~n=1,\ldots,N_m, ~r=1,\ldots,R\},
$$
 all with Hurst index $\widehat{H}_{m}^{(1)}$.  Obtain
 \[
 \mbox{\it
  the estimated Hurst index $\widehat{h}_{mr}^{(1)}$ of $\{G_{mr}(n),~n=1,\ldots,N_m\}$
}\]
and
  \[
  \mbox{\it the estimated Hurst index $\widehat{h}_{mr}^{(2)}$ of $\{G_{mr}(n)^2,~n=1,\ldots,N_m\}$,
}\]
 for each $r=1,\ldots,R$.
 
 \smallskip
\item For each $m=1,\ldots,M$, compute
\[\delta_m=\widehat{H}_{m}^{(2)}-\max(\frac{1}{2},2\widehat{H}_{m}^{(1)}-1)\]
from the data, and compute
 $$
 \left\{\delta_{mr}:=\widehat{h}_{mr}^{(2)}-\max\left(\frac{1}{2},  2\widehat{h}_{mr}^{(1)}-1\right), ~ r=1,\ldots,R \right\}
 $$
  from the simulated series.
Then compute the relative number of times (percentile) that $\delta_{mr}$ is less than or equal to $\delta_m$ for $r=1,\ldots,R$, namely,
\[
P_m=\frac{1}{R} \#\{\delta_{mr}:~ \delta_{mr}\le \delta_m,~r=1,\ldots,R \}=\widehat{F}_{m,R}(\delta_m),
\]
where $\widehat{F}_{m,R}$ is the empirical CDF of  $\{\delta_{mr}:~r=1,\ldots,R\}$.

\smallskip
\item Construct the following contrast group: for each $m=1,\ldots,M$,    simulate a fractional Gaussian noise  sequence  
$$
\{X_m^*(n),~n=1,\ldots,N_m^*\}
$$
 with Hurst index  randomly sampled from $\{\widehat{H}_{m}^{(1)}, ~m=1,\ldots,M\}$, and length $N_m^*$ randomly sampled from $\{N_m,~m=1,\ldots,M\}$.
 
\smallskip
Then perform the preceding steps 1-4 replacing $\{X_m(n)\}$ by $\{X_m^*(n)\}$, from which one gets $\widehat{H}_{m}^{(1)*}$, $\widehat{H}_{m}^{(2)*}$, $\delta_m^*$ and $P_m^*$ that correspond to $\widehat{H}_{m}^{(1)}$, $\widehat{H}_{m}^{(2)}$, $\delta_m$ and  $P_m$ respectively.

In our study we set $R=200$.
\end{enumerate}
}

\medskip
\begin{Rem}\label{Rem:explain desgin}
We explain here the preceding study design.
Recall that $\widehat{H}_m^{(1)}$ is the Hurst index estimate of the time series $\{X_m(n)\}$ and $\widehat{H}_m^{(2)}$ is the Hurst index estimate of the squared time series $\{X_m(n)^2\}$.
As mentioned before, the goal is to examine whether  $\widehat{H}_{m}^{(1)}$ and $\widehat{H}_{m}^{(2)}$ behave according to (\ref{eq:H_G rank 2}). If they behave perfectly according to (\ref{eq:H_G rank 2}),
then
\[
\delta_m=\widehat{H}_m^{(2)}-\max\left(\frac{1}{2},2\widehat{H}_m^{(1)}-1\right)
\]
should be  zero. Both  $\widehat{H}_{m}^{(1)}$ and $\widehat{H}_{m}^{(2)}$  are random and thus  fluctuate as $m$ varies. To get a reference point, we  introduce a statistical contrast in Step 3, whereby we simulate $R$ fractional Gaussian noises series and measure their Hurst indices $\widehat{h}_{mr}^{(1)}$ and $\widehat{h}_{mr}^{(2)}$: the first index is for fractional Gaussian noise and the second is for its square. Since these Hurst indices are obtained from   fractional Gaussian noises, they indeed obey (\ref{eq:H_G rank 2}).  We want to see how    $\delta_m$, which is measured from data, compares to  the $\delta_{mr}=\widehat{h}_{mr}^{(2)}-\max\left(\frac{1}{2},  2\widehat{h}_{mr}^{(1)}-1\right)$ corresponding to fractional Gaussian noise. This leads us to focus on   $P_m$ instead of  $\delta_m$. One may view $P_m$ as a  ``standardized'' version of $\delta_m$ with respect to the contrast distribution constructed from fractional Gaussian noise, which makes comparison across different data items (different $m$)  more consistent.
More technical explanations are given below.

 Let $F_m(x)$ be the CDF of the  random  $\delta_m$.  Then $F_m\left(\delta_m\right)$ follows exactly a \emph{uniform distribution} on $[0,1]$.  If $\{X_{m}(n)\}$ were indeed generated by fractional Gaussian noise  with true Hurst index $\widehat{H}_{m}^{(1)}$,  then the empirical CDF $\widehat{F}_{m,R}$ in Step 4 is a good approximation of $F_m$. Therefore,  if $\{X_m(n)\}$ obeys (\ref{eq:H_G rank 2}) as the fractional Gaussian noise does, and $\widehat{H}_{m}^{(1)}$ is a reasonable estimate, then $P_m=\widehat{F}_{m,R}(\delta_m)$   in Step 4  is expected to  follow  a \emph{uniform distribution} on $[0,1]$ approximately.    On the other hand, if the $\delta_m$ computed from the data makes the distribution of $P_m=\widehat{F}_{m,R}(\delta_m) $ skewed towards $1$, then this indicates that $\delta_m$ tends to be larger than $\delta_{mr}$.

To account for the potential bias due to the estimation of the Hurst index, in Step 5 we replace our original   data    $\{X_m(n)\}$ by a second contrast group $\{X_m^*(n)\}$ made up of fractional Gaussian noise sequences with similar lengths  and Hurst indices. After repeating the same procedure on this contrast group,  we  can then compare the distribution (histogram) of $\{P_m\}$ obtained from the original   data   with the distribution of $\{P_m^*\}$ obtained from the contrast group.

These designs may be regarded as  simulation-assisted statistical tests where the null hypothesis is the relation (\ref{eq:H_G rank 2}).
\end{Rem}

\medskip

\begin{figure}
\centering
\includegraphics[scale=0.5]{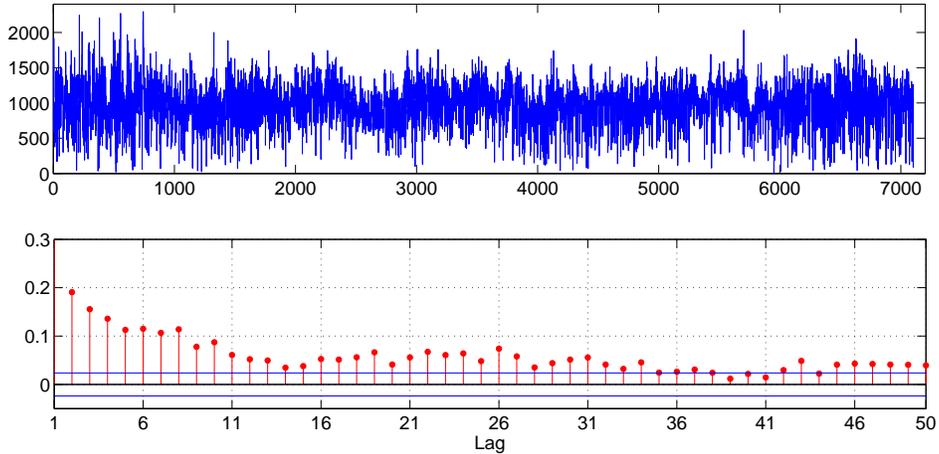}
\caption{\small{{\it Up:}\ Plot of the tree ring  time series extracted from ca506.crn in ITRDB. {\it Down:} autocorrelation plot. The variance aggregation estimate for the Hurst index of the data yields $\widehat{H}^{(1)}=0.7182$ and the Hurst index  for the centered and squared data yields $\widehat{H}^{(2)}=0.7217$; the local periodogram regression yields $\widehat{H}^{(1)}=0.7569$ and $\widehat{H}^{(2)}=0.7801$ respectively; the local Whittle estimate yields $\widehat{H}^{(1)}=0.7024$  and $\widehat{H}^{(2)}=0.7061$} respectively.} \label{fig:ca506}

\end{figure}
Now we describe the data we use. The \emph{tree ring width} in chronological order has been identified as one of the natural stationary time series data sets which exhibit long memory (see \citet{mandelbrot:wallis:1969:some} and \citet{pelletier:turcotte:1997:long}).
Since the tree ring width is largely affected by environmental factors, which is explored in dendrochronology  (see \citet{schweingruber:1996:tree}), it also reflects the   long-memory stationary fluctuation of the  ecological systems. We shall use the data compiled by The International Tree-Ring Data Bank (ITRDB, \url{ftp://ftp.ncdc.noaa.gov/pub/data/paleo/treering/chronologies/}) collected from Africa, Asia, Australia, Canada, Europe, Mexico, South America and USA, stored in the  Standard Chronology File (*.crn) format. For example, Figure \ref{fig:ca506} displays the time series extracted from the file ca506.crn in the data bank and its autocorrelation plot.
We further select the data according to the following criteria:
\begin{enumerate}[Cr{i}ter{i}on 1]
\item  The length of the time series is at least 300.
\item  The time series data is importable by the Tree-Ring Matlab Toolbox\footnote{\url{http://www.ltrr.arizona.edu/~dmeko/toolbox.html}} (data is usually importable if there is no missing value).
\item  The estimated Hurst index $\widehat{H}_{m}^{(1)}$ lies within the interval $[0.6,0.9]$\footnote{Ideally we want the selected data to be  stationary and long-range dependent. When the estimate is close to $0.5$, the data is likely to have short memory; when the estimate is close to $1$, it is likely to be  non-stationary.}.
\end{enumerate}
To be consistent, we also apply  Criterion 1 and Criterion 3 for to the contrast group $\{X_m^*(n),~n=1,\ldots,N_m\}$.

\medskip

We shall use the following three popular estimators of Hurst index:
\begin{itemize}
\item Variance aggregation estimator;
\item Local periodogram regression estimator (also known as GPH estimator);
\item Local Whittle estimator.
\end{itemize}
  For a description and empirical study of these estimators,  see \citet{taqqu:1995:estimators}. There are more sophisticated estimators, for example, the wavelet-type estimators (see, e.g., \citet{fay:2009:estimators}). To minimize finite-sample bias, these methods typically involve complicated choice of some tuning parameters. Since our study design has taken into account the potential bias of the estimator, we shall stick to the  three more elementary estimators aforementioned.
For the  variance aggregation estimator and the local periodogram regression (GPH estimator), we use the implementation by  Chu Chen (\url{http://www.mathworks.com/matlabcentral/fileexchange/19148-hurst-parameter-estimate}, and we use the default parameter settings);
For the local Whittle estimate, we use the implementation by Katsumi Shimotsu (
\url{http://shimotsu.web.fc2.com/Site/Matlab_Codes.html}), in which case we   choose the frequency cutoff threshold to be $[N^{2/3}]$ with $N$ being the length of the time series).

\medskip

\begin{figure}
\centering
\includegraphics[scale=0.4]{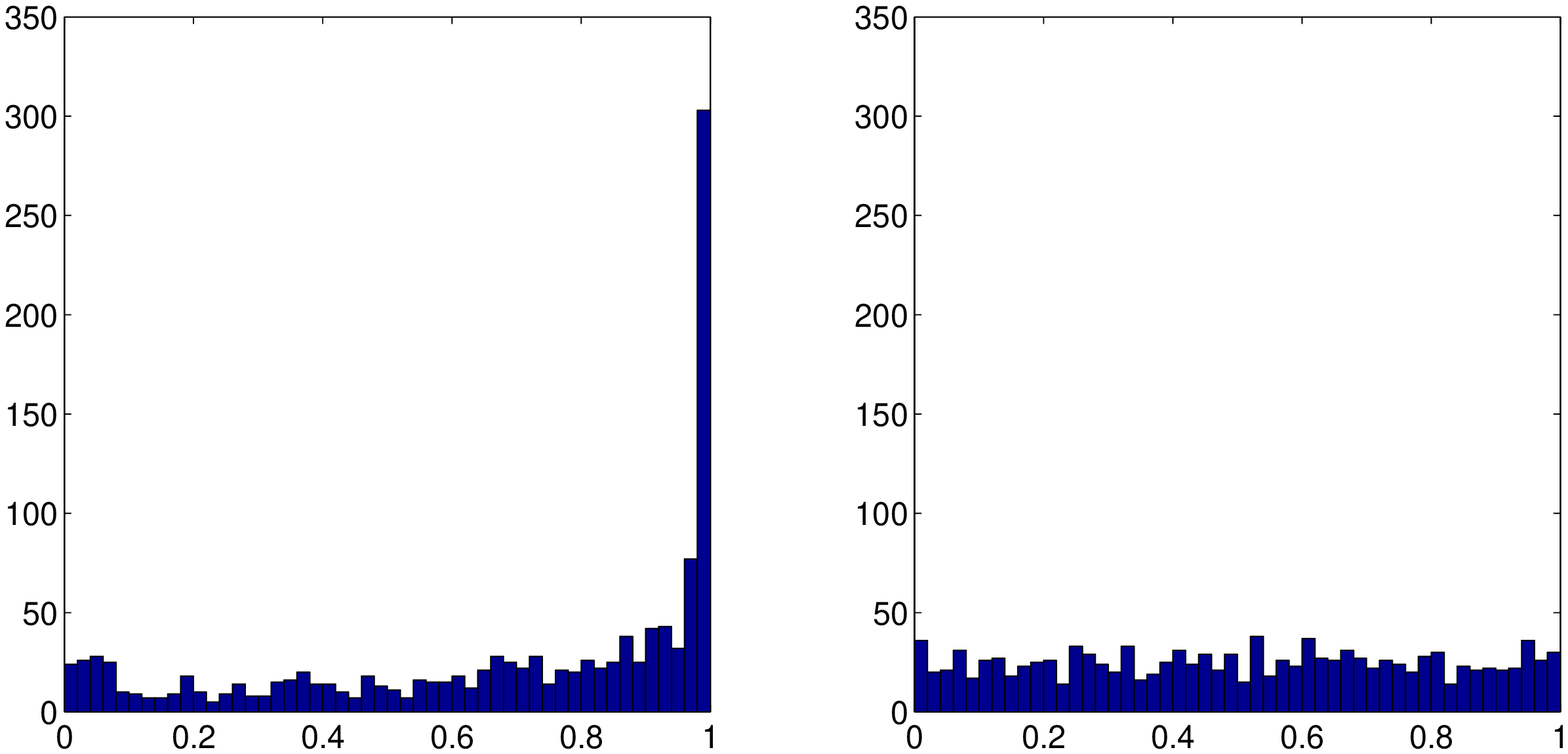}

\caption{\small{Histogram of $\{P_m\}$ (left) v.s.\ $\{P_m^*\}$ (right) from the Variance Aggregation Estimator}}\label{fig:aggvar_hist}


\centering
\includegraphics[scale=0.4]{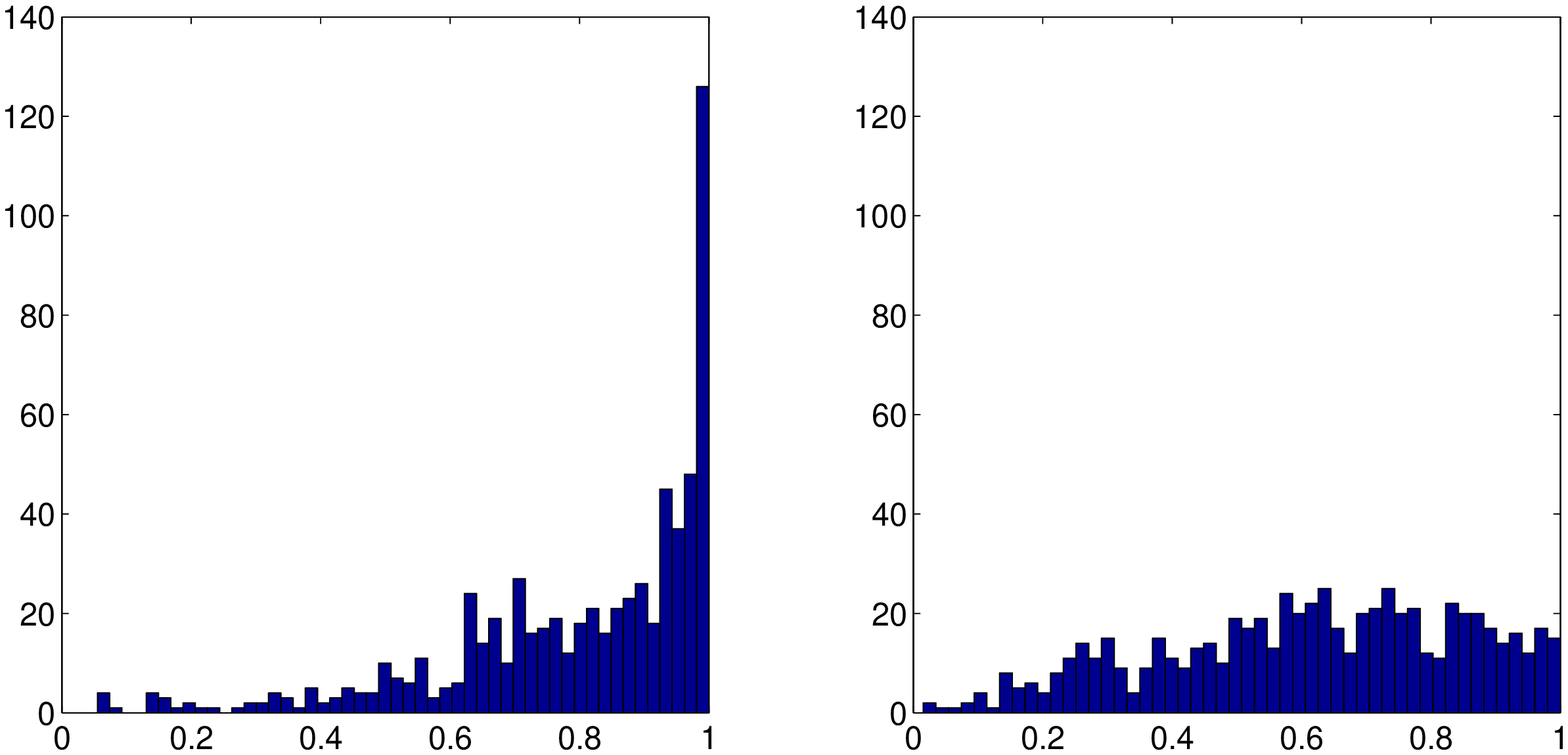}
\caption{\small{Histogram of $\{P_m\}$ (left) v.s.\ $\{P_m^*\}$ (right) from the Local Periodogram Regression}} \label{fig:per_hist}
%

\centering
\includegraphics[scale=0.4]{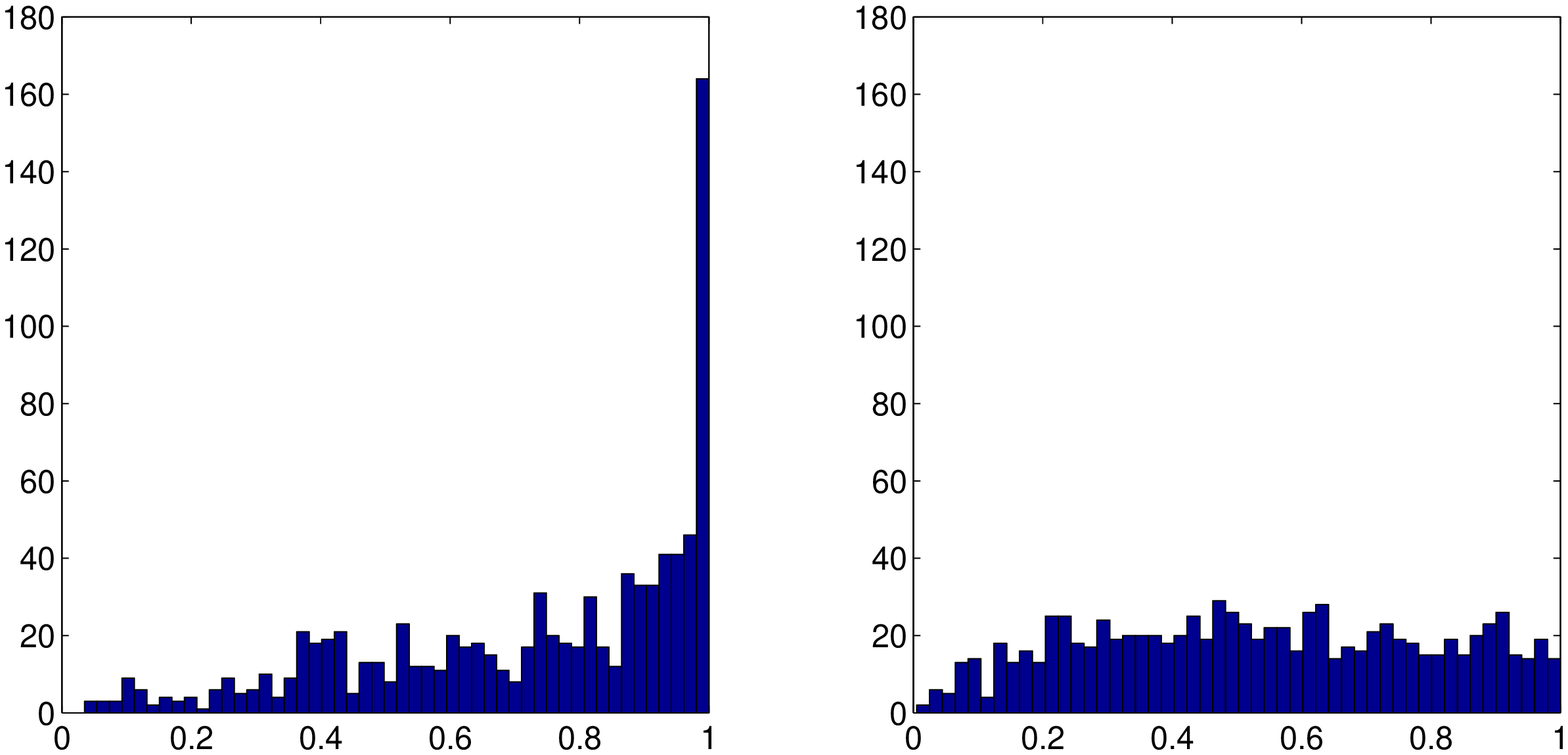}

\caption{\small{Histogram of $\{P_m\}$ (left) v.s.\ $\{P_m^*\}$ (right) from the Local Whittle Estimator}}\label{fig:whittle_hist}

\end{figure}

\noindent\textbf{Observations}:

\medskip

The graphs in the right-hand side of Figure \ref{fig:aggvar_hist}, \ref{fig:per_hist} and \ref{fig:whittle_hist} are as expected, namely, corresponding roughly to a uniform  distribution. This indicates that the procedure described in the study is reasonable.  In fact, the median of $P^*_m$ is roughly $50\%$ as it should be (see Table \ref{tab:summary}). As mentioned below, there may be a small bias when using the Local Periodogram Regression method (Figure \ref{fig:per_hist} (right)).   See also \citet{taqqu:teverovsky:1997:robustness} for an empirical discussion of  Whittle-type estimators.

Table \ref{tab:summary} summarizes some key statistics of the analysis  based on the three different estimators. One can see that for all three estimators, the median of  $\delta_m$ is consistently smaller than that of the contrast  $\delta_m^*$.
The median of  $P_m$ is  significantly  smaller than that of the contrast $P_m^*$.
Figure \ref{fig:aggvar_hist}, \ref{fig:per_hist} and \ref{fig:whittle_hist} plot the histograms of  $\{P_m\}$ and $\{P_m^*\}$  obtained via the three different estimators. Their results are similar: while $\{P_m^*\}$ are roughly uniformly distributed as expected, the histogram of $\{P_m\}$ is severely skewed towards $1$. The contrast in the skewness shows that the  $\delta_m$ computed from the tree ring data tends to be  larger than the $\{\delta_{mr}\}$ computed from the  fractional Gaussian noise.   In other words,  in the case of tree ring data, the Hurst index does \emph{not} tend to decrease  as much  after squaring  as  the case of  fractional Gaussian noise.

As mentioned in Remark \ref{Rem:explain desgin}, if the Hurst index estimate is unbiased, $P_m^*$ is expected to approximately follow a uniform distribution on $[0,1]$, so that the median is close to $1/2$. However, the estimation bias of Hurst index could distort this uniformity.  Indeed, in the Local Periodogram Regression case, the median of $P_m^*$ is 63.5\%.  But this is still in sharp contrast with the corresponding median of $P_m$ which is 86.25\% and hence significantly larger.
This  indicates that the data is not behaving like fractional Gaussian noise. Thus our design is effective despite the bias inherent in the estimation method.

\begin{table}
\begin{center}
\begin{footnotesize}
\begin{tabular}{|c|c|c|c|c|c|}
\hline
Estimator & Selected  number $M$ & Median   $\delta_m$  & Median   $\delta_m^*$ & Median $P_m$ & Median $P_m^*$ \\
\hline
Variance Aggregation & 1250 &  0.0786 &  0.0104 &  80.50\% &   51.00\% \\
\hline
Local Periodogram Regression & 658 & 0.0921 & -0.0204 &  86.25\% & 63.50\% \\
\hline
Local Whittle & 908 & 0.0496 & -0.0162 &  80.50\% & 52.50\% \\
\hline
\end{tabular}
\end{footnotesize}
\end{center}
\caption{Analysis Summary}\label{tab:summary}
\end{table}

\begin{figure}
\centering
\includegraphics[height=7cm,width=7cm]{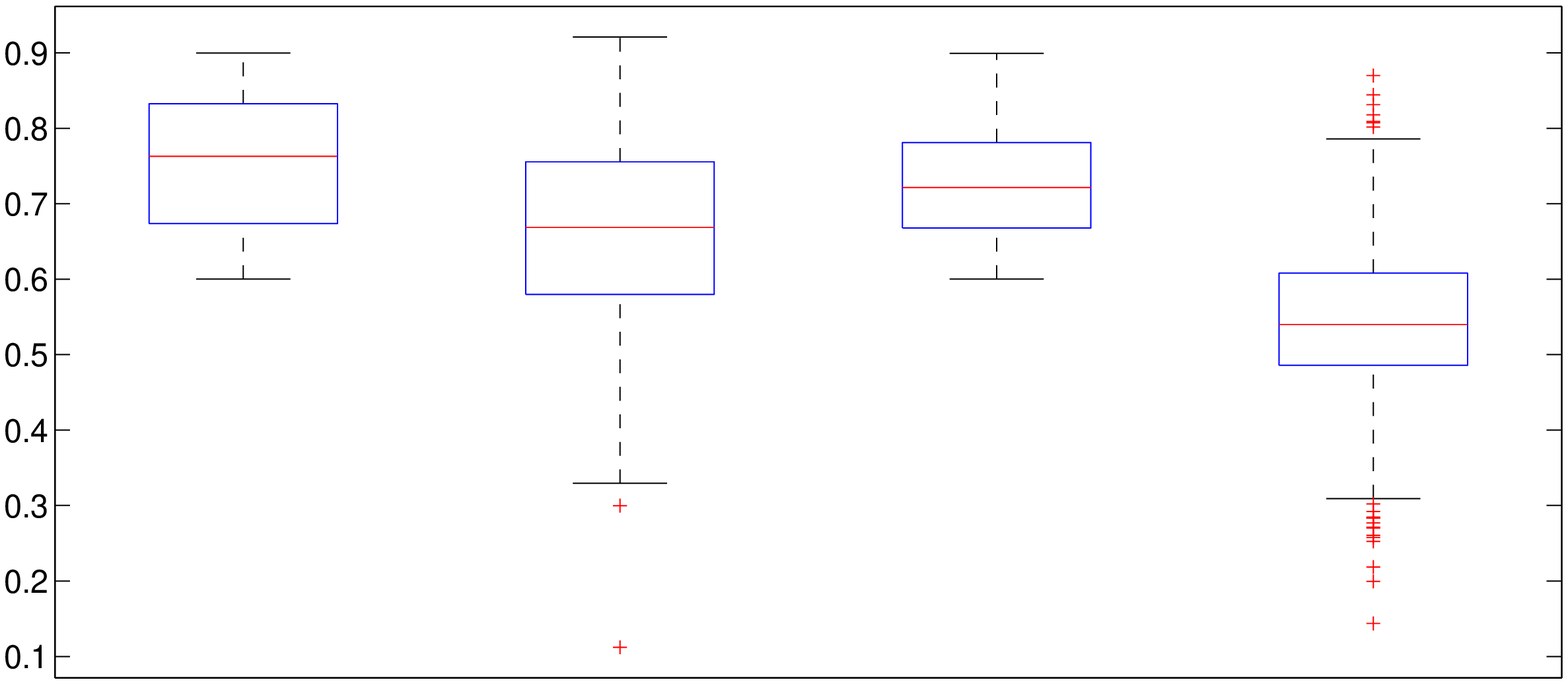}\\
\includegraphics[height=7cm,width=7cm]{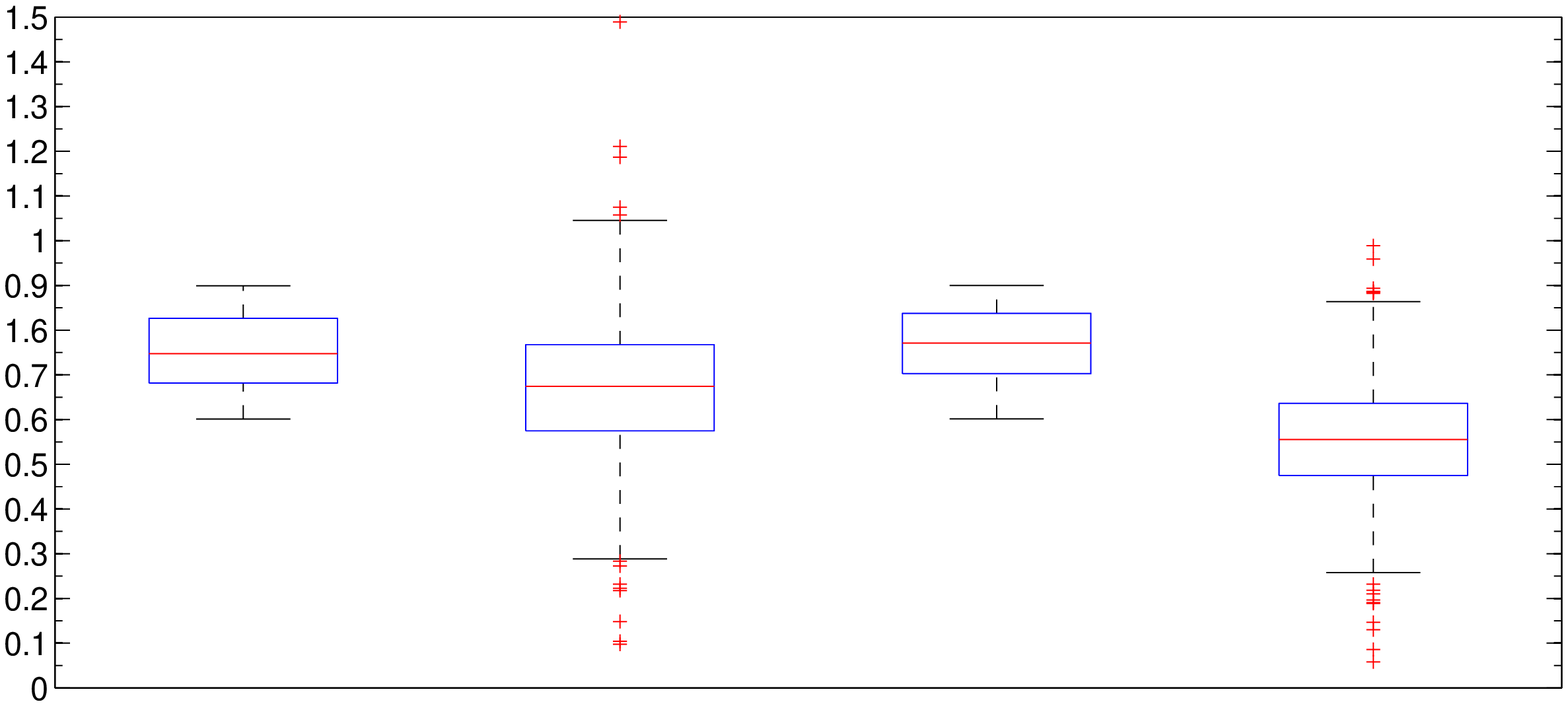}\\
\includegraphics[height=7cm,width=7cm]{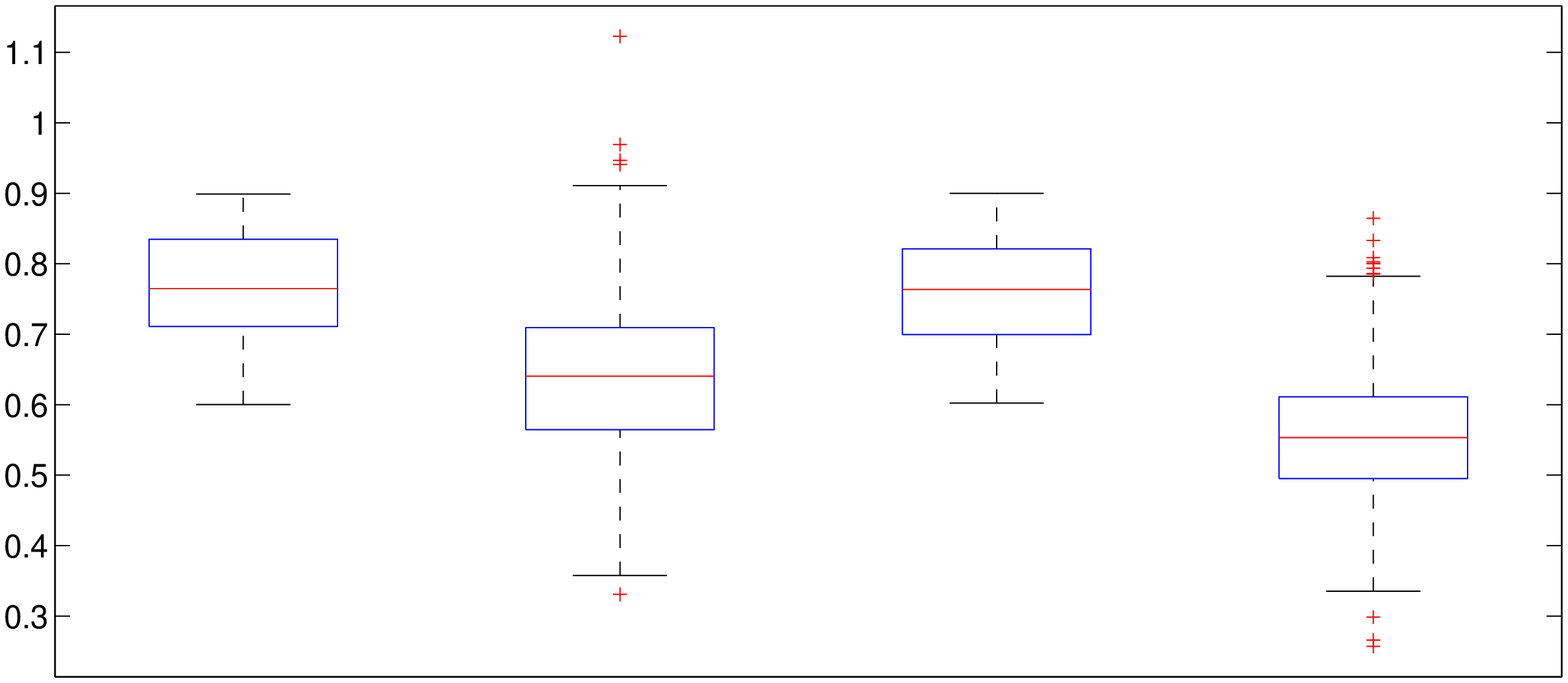}
\caption{\small{Top to bottom: variance aggregation estimator, local periodogram regression and local Whittle estimator.
In each boxplot, from left to right: $\widehat{H}_{m}^{(1)}$, $\widehat{H}_{m}^{(2)}$, $\widehat{H}_{m}^{(1)*}$ and $\widehat{H}_{m}^{(2)*}$.}}\label{fig:diff}
\end{figure}

\begin{Rem}\label{Rem:decay H2}
From the analysis above, we  conclude that relation (\ref{eq:H_G rank 2}), or more generally (\ref{eq:H_G general}),  may not make good prediction on real-life data. We note, however,  that the estimated Hurst index $\widehat{H}_{m}^{(2)}$ of $\{X_m(n)^2\}$  tends to be somewhat smaller than the estimated Hurst index $\widehat{H}_{m}^{(1)}$ of $\{X_m(n)\}$, although for the contrast group $\{X^*(n)\}$ the decrease from $\widehat{H}_{m}^{(1)*}$ to $\widehat{H}_{m}^{(2)*}$ is more significant. See Figure \ref{fig:diff}.
A possible explanation  is that
 although  $\{X_m(n)^2\}$  actually possesses rank $1$ and thus has the same Hurst index as $\{X_m(n)\}$,  many of the $\{X_m(n)\}$   may be close to a Gaussian (or linear) process. So they tend to exhibit somewhat the relation (\ref{eq:H_G rank 2}) when the sample size is moderate. See \citet{bai:taqqu:2017:hermite} for  an analysis of the interplay between the rank instability effect and the sample size.
\end{Rem}
\begin{Rem}
As a reviewer pointed out,    another explanation of the observations found in the study is that the data originally follows a model with a rank higher than $1$, in which case squaring  does not necessarily lead to  a higher-order rank.  Although this explanation is allowable in theory, it is less natural than the instability explanation.   The reviewer's  explanation relies on assuming a special model: the transformation of a Gaussian or linear process with higher-order rank, while ours indicates that a slight perturbation makes the formula (\ref{eq:H_G rank 2}) unrealistic in practice.
\end{Rem}

%
%

\section{Stability of limit theorems under weak dependence}\label{Sec:clt stable}
In this section, we demonstrate that the instability phenomenon appearing in the   limit theorems under long memory  does not typically occur in the short-memory case. This is  important  because it shows that the  transformation considered as ``perturbation'' in the previous section   usually does not make any qualitative difference in   short-memory  situations and hence may be safely negligible in large sample inference.

There are many ways to mathematically characterize weak dependence.
For an introduction to various notions of weak dependence of stationary processes and corresponding limit theorems,  we refer  to  \citet{doukhan:2003:models}. In this section, we shall mainly look at  the following three as examples:
\begin{enumerate}[(1)]
\item  Fast-decaying mixing coefficients under strong mixing conditions;
\item  Fast decaying covariance function
 in Gaussian subordination model (Theorem \ref{Thm:clt gauss});  \item Fast decaying physical dependence measure of \citet{wu:2005:nonlinear} in Bernoulli shift models.
\end{enumerate}
 The first is by far the most widely-used  notion for weak dependence which applies to very general stationary processes. The second is mentioned due to its close connection to the considerations in Section \ref{Sec:inst trans}.  The third is a convenient criterion under the Bernoulli shift framework which covers a wide range of concrete statistical models.

\smallskip
\subsection{Strong mixing conditions}
Suppose that $\{Y(n)\}$ is a stationary process with $\E [Y(n)]=0$ and $\Var[ Y(n)]=1$. Define the  $\sigma$-field $\mathcal{F}_{a}^b=\sigma\{Y(n): a\le  n\le b  \}$, where $-\infty\le  a\le b \le +\infty$. Given two  $\sigma$-fields $\mathcal{A}, \mathcal{B}$, one can define the following measure of dependence
\begin{equation}\label{eq:alpha dep measure}
\alpha(\mathcal{A},\mathcal{B})=\sup\{|P(A\cap B)-P(A)P(B)|:~ A\in \mathcal{A}, B\in \mathcal{B}\}.
\end{equation}
Then the $\alpha$-mixing coefficient of $\{X(n\})$, first introduced in \citet{rosenblatt:1956:central}, is defined as
\begin{equation*}
\alpha_Y(n)=\alpha\left( \mathcal{F}_{-\infty}^0, \mathcal{F}_{n}^\infty \right).
\end{equation*}
When $\alpha_Y(n)\rightarrow 0$ as $n\rightarrow\infty$, we say that $\{Y(n)\}$ is  \emph{strong mixing}. If one assumes that $\alpha_Y(n)$ decays to zero fast enough together with some other regularity conditions, then a central limit theorem for   $X(n)$ can be established. We state, as an example, the following  central limit theorem due to \citet{ibragimov:1962:some} and \citet{herrndorf:1984:functional}.
\begin{Thm}\label{Thm:strong mixing clt}
If $\E |Y(n)|^{2+\delta}<\infty$ for some $\delta>0$ and
\begin{equation}\label{eq:alpha summable}
\sum_{n=1}^\infty \alpha_Y(n)^{\delta/(2+\delta)}<\infty,
\end{equation}
then
\begin{equation*}
\frac{1}{\sqrt{N}} \sum_{n=1}^{[Nt]} \Big( Y(n)-\E Y(n) \Big) \Rightarrow \sigma B(t).
\end{equation*}
where $B(t)$ is a standard Brownian motion, $\Rightarrow $ stands for weak convergence in  $D[0,1]$, and
\[
\sigma^2=\sum_{n=-\infty}^\infty \Cov[Y(n),Y(0)].
\]
\end{Thm}

Now consider the   transformation    \[
X(n)=F(Y(n),\ldots,Y(n-l)).\] Let us compare $\alpha_X$ and $\alpha_Y$.
Since $X(n)\in \mathcal{F}_{n-l}^n$, it is easily deduced that for $n>l$, the $\alpha$-mixing coefficient of $\{X(n)\}$ satisfies
\begin{equation}\label{eq:alpha y bound}
\alpha_X(n) \le \alpha_Y(n-l).
\end{equation}
The relation (\ref{eq:alpha y bound}) means that the dependence measured by the $\alpha$-mixing coefficient after the perturbing tranform $F(\cdot)$ cannot exceed that of the original process $Y(n)$ (up to a fixed lag $l$). In particular, relation (\ref{eq:alpha summable}) holds for $\alpha_X(n).$ One then only needs  $\E |X(n)|^{2+\delta}<\infty$ (which is the case if $F(\cdot)$ has at most linear growth) for Theorem \ref{Thm:strong mixing clt}  to hold.

There are different mixing coefficients than $\alpha(n)$, obtained by modifying the measure of dependence between the $\sigma$-fields in (\ref{eq:alpha dep measure}),
for example, the $\phi$-mixing coefficient defined through
\[
\phi(\mathcal{A},\mathcal{B})= \sup\{|P(A|B)-P(A)|:~ A\in \mathcal{A}, B\in \mathcal{B}, P(B)>0\},
\]
the $\rho$-mixing coefficient defined through
\[
\rho(\mathcal{A},\mathcal{B})= \sup\{\mathrm{Corr}(X,Y):~ X\in L^2(\mathcal{A}), Y\in L^2(\mathcal{B})\},
\]
and so on. In general, as long as a dependence measure $m(\cdot,\cdot)$ is non-increasing with respect to set inclusion and the mixing coefficient is defined as $m(n)=m(\mathcal{F}_{-\infty}^0,\mathcal{F}_n^\infty)$, then
 a relation as (\ref{eq:alpha y bound}) always  holds.

Hence, the central limit theorems under strong mixing conditions is robust against a  transformation perturbation.

\smallskip
\subsection{Gaussian  subordination}
Let $\{Y(n)\}$ be a stationary \emph{Gaussian} process, and let
\[
X(n)=F(Y(n),\ldots,Y(n-l)).
\]
When the covariance function of $Y(n)$ decays fast enough, a central limit theorem always holds for $X(n)$.
In particular, we have the following result which is a consequence  of \citet{ho:sun:1987:central}.
\begin{Thm}\label{Thm:clt gauss}
Suppose that $\E X(n)^2<\infty$ and
\begin{equation}\label{eq:summable cov X}
\sum_{n=-\infty}^\infty \left|\Cov\left[Y(n),Y(0)\right] \right|<\infty.
\end{equation}
Then one has
\begin{equation*}\label{eq:clt gauss}
\frac{1}{\sqrt{N}} \sum_{n=1}^{[Nt]} \Big(X(n)-\E X(n) \Big)\ConvFDD \sigma B(t),
\end{equation*}
where $B(t)$ is a standard Brownian motion and
\[
\sigma^2=\sum_{n=-\infty}^\infty \Cov[X(n),X(0)].
\]
\end{Thm}

Theorem \ref{Thm:clt gauss} directly expresses the robustness of the central limit theorem against  transformation perturbation when the short memory condition (\ref{eq:summable cov X}) is imposed on $Y(n)$.

\smallskip
\subsection{Bernoulli shift}
Let $\{\epsilon_i\}$ be an i.i.d. sequence of random variables with mean $0$ and variance $1$.
Consider the Bernoulli shift model
\begin{equation}\label{eq:bernoulli shift}
Y(n)=G_Y(\epsilon_n,\epsilon_{n-1},\ldots),
\end{equation}
where $G_Y$ is a non-random measurable function.
 This specification covers not only the causal linear process (\ref{eq:linear}), but also many nonlinear time series models obtained as  solutions of difference equations involving $\epsilon_i$.

\citet{wu:2005:nonlinear} introduced the following so-called \emph{physical dependence measure} for a process $\{Y(n)\}$ specified by (\ref{eq:bernoulli shift}). Let $\epsilon_0^*$ be a random variable independent of $\{\epsilon_i\}$ and having the same distribution as $\epsilon_0$. Define
\begin{equation}\label{eq:delta_2}
\delta_2^X(n)=\|G_Y(\epsilon_n,\ldots,\epsilon_1,\epsilon_0,\epsilon_{-1},\ldots) -G_Y(\epsilon_n,\ldots,\epsilon_1,\epsilon_0^*,\epsilon_{-1},\ldots)\|_{L^2(\Omega)}.
\end{equation}
If (\ref{eq:bernoulli shift}) is interpreted as a nonlinear system with input $\{\epsilon_n\}$ and output $\{Y(n)\}$,  then $\delta_2^Y(n)$ in (\ref{eq:delta_2})  measures the influence of the lag-$n$ input $\epsilon_0$  on the current output $Y(n)$.

With $\delta_2^Y(n)$, one can state the following central limit theorem, which is a consequence of Theorem 1 and 3 of  \citet{wu:2005:nonlinear}.
\begin{Thm}\label{Thm:clt Wu}
Suppose that
\begin{align}\label{eq:weak dep Wu}
\sum_{n=1}^\infty\delta_2^X(n)<\infty.
\end{align}
Then one has
\begin{equation*}
\frac{1}{\sqrt{N}} \sum_{n=1}^{[Nt]} \Big( X(n)-\E X(n)\Big) \ConvFDD \sigma B(t),\quad t\ge 0,
\end{equation*}
where $B(t)$ is a standard Brownian motion,   and
\[
\sigma^2=\sum_{n=-\infty}^\infty \Cov[X(n),X(0)].
\]
\end{Thm}
\begin{Rem}
The criterion (\ref{eq:weak dep Wu}) is typically  easier to check for a specific Bernoulli shift model than the criteria based on strong mixing conditions (see Theorem \ref{Thm:strong mixing clt}), while still providing numerous statistical applications.
\end{Rem}

Now we consider the   transformation perturbation.
Let
\[
X(n)=F(Y(n),\ldots,Y(n-l+1))=:G_X(\epsilon_{n},\epsilon_{n-1},\ldots).
\]
We need to assume some smoothness condition (compare with the arguments of Claim \ref{Cla:Herm}) on the  perturbation function $F(x_1,\ldots,x_l)$. In particular, suppose that $F(\cdot)$ is Lipschitz, that is,
\begin{equation}\label{eq:lipschitz}
|F(x_1,\ldots,x_l) -F(y_1,\ldots,y_l)| \le  C_F\sum_{i=1}^l |x_i-y_i|
\end{equation}
for some constant $C_F\ge 0$.

Setting $\mathbd{\epsilon}_n=(\epsilon_n,\ldots,\epsilon_1,\epsilon_0,\epsilon_{-1},\ldots)$ and $\mathbd{\epsilon}_n^*=(\epsilon_n,\ldots,\epsilon_1,\epsilon_0^*,\epsilon_{-1},\ldots)$, one has   by (\ref{eq:lipschitz}) that
\begin{align*}
\left|G_X(\mathbd{\epsilon}_n) -G_X(\mathbd{\epsilon}_n^*)\right|
\le C_F\sum_{i=0}^{l-1} |G_Y(\mathbd{\epsilon}_{n-i})-G_Y(\mathbd{\epsilon}_{n-i}^*) |.
\end{align*}
Therefore, if  $\delta_2^X(n)$  and $\delta_2^Y(n)$ are the physical dependence measures  of $\{X(n)\}$ and $\{Y(n)\}$ respectively, then
\begin{align*}
\delta_2^X(n)=\|G_X(\mathbd{\epsilon}_n) -G_X(\mathbd{\epsilon}_n^*)\|_{L^2(\Omega)}\le C_F \sum_{i=0}^{l-1} \delta_2^Y(n-i).
\end{align*}
Hence if $\{Y(n)\}$ satisfies the short memory condition 
$$
\sum_{n=1}^\infty\delta_2^Y(n)<\infty,
$$
 then so  does $\{X(n)\}$.
 This shows the robustness of Theorem \ref{Thm:clt Wu} against a perturbation by any Lipschitz   transformation.
\begin{Rem}
The proof of Theorem \ref{Thm:clt Wu} is based on a martingale difference approximation method and resorts to the martingale difference central limit theorem. We note, however, that the martingale difference central limit theorem  is itself not robust against  transformation, since the martingale difference structure  in general can be easily disturbed by a  transformation. For example, in the stochastic volatility-type models, e.g., the LARCH($\infty$) model (\citet{giraitis:2004:larch}), the return sequence $X(n)$ is a martingale difference, while $|X(n)|$ can exhibit long memory (see \citet{beran:2013:long}, Chapter 4.2.8.).
\end{Rem}

\begin{Rem}
Using similar arguments, one can show that the $\theta$-weak dependence criterion (whose definition involves bounded Lipschitz transformation) introduced by \citet{doukhan:louhichi:1999:new},  enjoys a robustness against bounded Lipschitz transformations.
\end{Rem}

\section{Conclusion and suggestions}\label{Sec:conclusion}

In this paper,   we discussed the instability  issue of Hermite rank and other related ranks   appearing in   limit theorems under long memory. We argued that a  rank  greater than $1$ can be  disturbed  by a transformation and only  a rank equal to $1$ is stable.   We provided empirical evidence supporting this argument. Such an instability feature has important   statistical implications. In particular, assuming a higher-order   rank when it is really not there may result in underestimating  the order of fluctuation of the statistic of interest.

To address this issue   we briefly indicate here some   suggestions for performing valid inference.
As  illustrated, particularly in Section \ref{Sec:examples}, one may adopt   the assumption that the rank is always $1$, regardless of any  nonlinear transformation resulting from the statistical procedure. Here the rank should be understood in a generalized sense, taking into account   situations as (\ref{eq:rank whittle}).
 Some studies have implicitly done so, although without giving an explanation (see, e.g., \citet{beran:1991:m} and \citet{shao:2011:simple}). Recently \citet{beran:2016:testing} designed a statistical test based on resampling to distinguish Hermite rank $1$ and a higher-order Hermite in the model (\ref{eq:inst trans}).


Another appealing way out, is to redesign the statistical procedure in a way as to avoid using the fixed-rank limit theorems for inference directly. This may be achieved by combining  re-sampling
method (see, e.g, \citet{hall:1998:sampling}, \citet{nordman:lahiri:2005:validity}, \citet{zhang:2013:block}), \citet{bai:taqqu:2016:validity}), together with suitable self-normalization technique (see, e.g., \citet{shao:2010:self} and \citet{shao:2011:simple}). We refer the reader to \citet{jach:2012:subsampling} \citet{betken:Wendler:2015:subsampling} and \citet{bai:taqqu:zhang:2015:unified} for approaches of this type.

\appendix

\section{Non-instantaneous transformation of the Gaussian}\label{Sec:non-inst trans}
Let $\{Y(n)\}$  be   a standardized stationary long-memory  Gaussian process with Hurst index $H$. We  extend  here the discussion on instantaneous transformation  (\ref{eq:inst trans}) to the \emph{non-instantaneous} transformation
\begin{equation}\label{eq:non-inst tran}
X(n)=F\big(Y(n),Y(n-1),\ldots,Y(n-l)\big),
\end{equation}
where $X(n)\in L^2(\Omega)$ and $l$ is a finite positive integer. Since the non-instantaneous case  is much less treated in the literature, we shall    introduce  in this section   the relevant results in \citet{dobrushin:major:1979:non}, and show that the arguments developed in Section \ref{Sec:inst trans} continue to be valid.

It is well-known that the Gaussian $Y(n)$ admits the spectral representation (see, e.g., \citet{dobrushin:major:1979:non})
\begin{equation}\label{eq:spec rep}
Y(n)=\int_{(-\pi,\pi]} e^{inx} W_Y(dx),
\end{equation}
where $W_Y(dx)$ is a complex-valued Gaussian measure satisfying
\begin{equation}\label{eq:F_Y}
\E |W_Y(dx)|^2=F_Y(dx)
\end{equation}
and $F_Y(\cdot)$ is the spectral distribution\footnote{Do not confuse $F_Y$ in (\ref{eq:F_Y}) with $F$ in (\ref{eq:non-inst tran}).} of $Y(n)$.
Then $X(0)$ has following Wiener-It\^o expansion (see \citet{dobrushin:major:1979:non}, formula (6.1), or \citet{janson:1997:gaussian}, Theorem 7.61):
\begin{equation}\label{eq:wiener ito raw}
X(0)-\E X(0)= \sum_{m=1}^\infty \int_{(-\pi,\pi]^m}'' \alpha_m(x_1,\ldots,x_m) W_Y(dx_1)\ldots W_Y(dx_m),
\end{equation}
where the double prime $''$ indicates the exclusion of the hyper-diagonals $x_p=\pm x_q$ in the multiple stochastic integral. Here $\alpha_m(\cdot)$'s are a.e.\ unique complex-valued functions in satisfying
\[\alpha_m(x_1,\ldots,x_k)=\overline{\alpha_m(-x_1,\ldots,-x_m)},
\]  and
\[
\sum_{m=1}^\infty  m!\|\alpha_m\|_{L^2((-\pi,\pi]^m,F_Y^{\otimes m})}^2 <\infty,
\]
where
\[
\|\alpha_m\|_{L^2((-\pi,\pi]^m,F_Y^{\otimes m})^2} =\int_{(-\pi,\pi]^m} |\alpha_m(x_1,\ldots,x_m)|^2 F_Y(dx_1)\ldots F_Y(dx_m).
\]
The \emph{Hermite rank} of $X(n)$ (or say the Hermite rank of $F(\cdot)$ with respect to $\{Y(n)\}$) is defined as
\begin{equation}\label{eq:gen herm rank}
\inf\left\{m\ge 1: ~\|\alpha_m\|_{L^2((-\pi,\pi]^m,F_Y^{\otimes m})} \neq 0 \right\},
\end{equation}

The Hermite rank in (\ref{eq:gen herm rank})
is also equal to  (see \citet{dobrushin:major:1979:non} Remark 6.3)
\begin{equation}\label{eq:herm rank noninst explicit}
\inf\left\{m\ge 1:~\E\Big[ \big(X(0)-\E X(0)\big)Y(n)^{m} \Big] \neq 0 ~\text{ for some }n\in \mathbb{Z}\right\}.
\end{equation}
This should be compared to (\ref{eq:alter herm rank def}).

By Remark 6.1 of \citet{dobrushin:major:1979:non}, the   a.e.\ unique function $\alpha_m(\cdot)$ can further be chosen to be continuous, which we shall assume throughout below.
We are now ready to state  the following generalization of Theorem \ref{Thm:nclt gaussian}, which follows from \citet{dobrushin:major:1979:non} Theorem 3,  Remark 6.3 and Remark 6.4.
\begin{Thm}\label{Thm:nclt gen}
Suppose that   $X(n)=F(Y(n),\ldots,Y(n-l))$, and that the Hermite rank in the sense of (\ref{eq:gen herm rank})  is $k$, and that the Hurst index $H$ of $\{Y(n)\}$  satisfies
\[
H>1-\frac{1}{2k}.
\]
Suppose also that $\alpha_k(\cdot)$ in (\ref{eq:wiener ito raw}) satisfies
\begin{equation}\label{eq:coef nonzero}
\alpha_k(0,\ldots,0)\neq 0.
\end{equation}
Then $\{X(n)\}$   has long memory with Hurst index:
\[H_F=(H-1)k+1\in \left(\frac{1}{2},1\right).
\]
Furthermore, as $N\rightarrow\infty$, we have
\begin{equation}\label{eq:Herm limit non-inst}
\frac{1}{N^{H_F}}\sum_{n=1}^{[Nt]}\Big( X(n)-\E X(n) \Big) \Rightarrow c \alpha_k(0,\ldots,0) Z_{H_F,k}(t),
\end{equation}
for some $c\neq 0$, where $Z_{H_F,k}(t)$ is the Hermite process  in (\ref{eq:Herm process}).
\end{Thm}
\begin{Rem}
In contrast to Theorem \ref{Thm:nclt gaussian} where the constant $c$ in (\ref{eq:Herm limit}) is always nonzero, in the non-instantaneous case we need to assume in addition the condition (\ref{eq:coef nonzero}).  If $\alpha_k(0,\ldots,0)=0$, then  (\ref{eq:Herm limit non-inst}) tells nothing more than that the normalization $N^{-H_F}$ is too strong. In this case, terms with order greater than $k$ may contribute to the asymptotic distribution as well. For example, if in (\ref{eq:non-inst tran})  we let
\[X(n)=H_1(Y(n))-H_1(Y(n-1))+H_2(Y(n))=Y(n)-Y(n-1)+Y(n)^2-1.\]
Using the spectral representation (\ref{eq:spec rep}) and  \citet{major:2014:multiple} Theorem 4.3,  we have
\[
X(0)=\int_{(-\pi,\pi]} ( 1-e^{-ix}  ) W_Y(dx)  +\int_{(-\pi,\pi]^2}''    W_Y(dx_1)W_Y(dx_2)
\]
so that $\alpha_1(x)=1-e^{-ix}$ and $\alpha_1(0)=0$.  On the  other hand, the Hermite rank of $X(n)$ is  $k=1$ in view of (\ref{eq:wiener ito raw}). Now
\[
\sum_{n=1}^N X(n) =Y(N)-Y(0)+ \sum_{n=1}^N H_2(Y(n)).
\]
Since $Y(n)$ is stationary, $N^{-H}[Y(N)-Y(0)]\ConvP 0$, and thus
only the term $\sum_{n=1}^N H_2(Y(n))$ contributes to in the limit.  Hence the limit of suitably normalized $\sum_{n=1}^N X(n)$  can  be either a Brownian motion if $H\le 3/4$ or a Hermite process of order $2$ if $H>3/4$ in view of Theorem \ref{Thm:nclt gaussian}.
\end{Rem}

\begin{Rem}
Now arguing  as in Section \ref{Sec:inst trans}, one notes that  a Hermite rank higher than $1$ in this non-instantaneous context is also unstable.
Recall that the role of $F(\cdot)$ in (\ref{eq:non-inst tran}), as in Section \ref{Sec:inst trans}, is to account for an uncontrollable perturbation of the Gaussian model.
 Suppose that $G(\cdot)$ is a function determined by the statistical procedure of interest. Then one can formulate a statement  parallel to Claim \ref{Cla:Herm}. So the part of Theorem \ref{Thm:nclt gen} which is most likely of statistical relevance is just the case $k=1$, where the limit is fractional Brownian motion  and the normalization is $N^{-H}$.  Note that this non-instantaneous consideration    includes   not only $G(X(n))$ with $X(n)$ defined in (\ref{eq:non-inst tran}), but also
 the case where $G(\cdot)$ is a finite-dimensional multivariate function of  the observed time series $\{X(n)\}$, for example $G(X(n),\ldots,X(n-p))=X(n)X(n-p)$, a term which appears in the sample covariance.

Arguing as in Claim \ref{Cla:Herm},  condition (\ref{eq:coef nonzero})   should be expected    to typically hold in practice.
\end{Rem}

\begin{Rem}
Using the full generality of Theorem 3 of \citet{dobrushin:major:1979:non},
it is even possible to consider the case $l=\infty$ in (\ref{eq:non-inst tran}), namely, including dependence on the infinite past.    In this case, however, one encounters major technical difficulties since $F(\cdot)$ with $l=\infty$ may alter the long memory property of $Y(n)$,  for example, if $F(\cdot)$ is a linear filter with a slow power-law decay (see, e.g., Section \ref{Sec:linear} below).  On the other hand, one may be satisfied with the restriction to $l<\infty$ since $F(\cdot)$  has been introduced  only to  account for a small  perturbation of the Gaussian model, in which case the argument of $F(\cdot)$  is  not expected to  stretch to the infinite past.
\end{Rem}

\begin{Rem}
This discussion can also be  extended to the case where $\mathbf{Y}(n)$ is a vector-valued Gaussian stationary noise and $\mathbf{X}(n)$ is also vector-valued. See, e.g., \citet{denaranjo:1993:non} and \citet{arcones:1994:limit}.
\end{Rem}
\begin{Rem}
We mention that   the extension of Theorem \ref{Thm:nclt linear} to non-instantaneous transformation of linear processes, that is, an analog of  Theorem \ref{Thm:nclt gen} when $X(n)$ is linear, is still open. Only  central limit theorems involving non-instantaneous filter of linear processes have been considered (see \citet{wu:2002:central} and \citet{cheng:Ho:2008:berry}).
\end{Rem}

\bigskip
\noindent \textbf{Acknowledgment.} We thank an associate editor and two referees for their insightful comments.  This work was partially supported by the NSF grant  DMS-1309009  at Boston University.

\bibliographystyle{plainnat}
\bibliography{Bib}

\begin{thebibliography}{86}
\providecommand{\natexlab}[1]{#1}
\providecommand{\url}[1]{\texttt{#1}}
\expandafter\ifx\csname urlstyle\endcsname\relax
  \providecommand{\doi}[1]{doi: #1}\else
  \providecommand{\doi}{doi: \begingroup \urlstyle{rm}\Url}\fi

\bibitem[Arcones(1994)]{arcones:1994:limit}
M.A. Arcones.
\newblock Limit theorems for nonlinear functionals of a stationary {G}aussian
  sequence of vectors.
\newblock \emph{The Annals of Probability}, pages 2242--2274, 1994.

\bibitem[Avram(1988)]{avram:1988:bilinear}
F.~Avram.
\newblock On bilinear forms in {G}aussian random variables and {T}oeplitz
  matrices.
\newblock \emph{Probability Theory and Related Fields}, 79\penalty0
  (1):\penalty0 37--45, 1988.

\bibitem[Avram and Taqqu(1987)]{avram:1987:noncentral}
F.~Avram and M.S. Taqqu.
\newblock Noncentral limit theorems and {A}ppell polynomials.
\newblock \emph{The Annals of Probability}, 15\penalty0 (2):\penalty0 767--775,
  1987.

\bibitem[Bai and Taqqu(2013)]{bai:taqqu:2013:1-multivariate}
S.~Bai and M.S. Taqqu.
\newblock Multivariate limit theorems in the context of long-range dependence.
\newblock \emph{Journal of Time Series Analysis}, 34\penalty0 (6):\penalty0
  717--743, 2013.

\bibitem[Bai and Taqqu(2016)]{bai:taqqu:2016:validity}
S.~Bai and M.S. Taqqu.
\newblock On the validity of resampling methods under long memory.
\newblock \emph{To appear in The Annals of Statistics}, 2016.

\bibitem[Bai and Taqqu(2017)]{bai:taqqu:2017:hermite}
S.~Bai and M.S. Taqqu.
\newblock Some properties of the {H}ermite rank.
\newblock \emph{Preprint, see ArXiv http://arxiv.org/abs/1710.01612}, 2017.

\bibitem[Bai et~al.(2016)Bai, Taqqu, and Zhang]{bai:taqqu:zhang:2015:unified}
S.~Bai, M.S. Taqqu, and T.~Zhang.
\newblock A unified approach to self-normalized block sampling.
\newblock \emph{Stochastic Processes and their Applications}, 126\penalty0
  (8):\penalty0 2465--2493, 2016.

\bibitem[Beran(1991)]{beran:1991:m}
J.~Beran.
\newblock M estimators of location for {G}aussian and related processes with
  slowly decaying serial correlations.
\newblock \emph{Journal of the American Statistical Association}, 86\penalty0
  (415):\penalty0 704--708, 1991.

\bibitem[Beran and Ghosh(1991)]{beran:gosh:1991:slowly}
J.~Beran and S.~Ghosh.
\newblock Slowly decaying correlations, testing normality, nuisance parameters.
\newblock \emph{Journal of the American Statistical Association}, 86\penalty0
  (415):\penalty0 785--791, 1991.

\bibitem[Beran and Weiersh{\"a}user(2011)]{beran:2011:spline}
J.~Beran and A.~Weiersh{\"a}user.
\newblock On spline regression under {G}aussian subordination with long memory.
\newblock \emph{Journal of Multivariate Analysis}, 102\penalty0 (2):\penalty0
  315--335, 2011.

\bibitem[Beran et~al.(2013)Beran, Feng, Ghosh, and Kulik]{beran:2013:long}
J.~Beran, Y.~Feng, S.~Ghosh, and R.~Kulik.
\newblock \emph{Long-Memory Processes: Probabilistic Properties and Statistical
  Methods}.
\newblock Springer, 2013.

\bibitem[Beran et~al.(2016)Beran, M{\"o}hrle, and Ghosh]{beran:2016:testing}
J.~Beran, S.~M{\"o}hrle, and S.~Ghosh.
\newblock Testing for {H}ermite rank in {G}aussian subordination processes.
\newblock \emph{Journal of Computational and Graphical Statistics}, 25\penalty0
  (3):\penalty0 917--934, 2016.

\bibitem[Betken and Wendler(2015)]{betken:Wendler:2015:subsampling}
A.~Betken and M.~Wendler.
\newblock Subsampling for general statistics under long range dependence.
\newblock \emph{arXiv preprint arXiv:1509.05720}, 2015.

\bibitem[Bingham et~al.(1989)Bingham, Goldie, and
  Teugels]{bingham:goldie:teugels:1989:regular}
N.H. Bingham, C.M. Goldie, and J.L. Teugels.
\newblock \emph{Regular Variation}.
\newblock Encyclopedia of Mathematics and Its Applications. Cambridge
  University Press, 1989.

\bibitem[Breuer and Major(1983)]{breuer:major:1983:central}
P.~Breuer and P.~Major.
\newblock Central limit theorems for non-linear functionals of {G}aussian
  fields.
\newblock \emph{Journal of Multivariate Analysis}, 13\penalty0 (3):\penalty0
  425--441, 1983.

\bibitem[Brockwell and Davis(1991)]{brockwell:1991:time}
P.J. Brockwell and R.A. Davis.
\newblock \emph{Time Series: Theory and Methods}.
\newblock Springer, 1991.

\bibitem[Chambers and Slud(1989)]{chambers:slud:1989:central}
D.~Chambers and E.~Slud.
\newblock Central limit theorems for nonlinear functionals of stationary
  {G}aussian processes.
\newblock \emph{Probability Theory and Related Fields}, 80\penalty0
  (3):\penalty0 323--346, 1989.

\bibitem[Cheng and Robinson(1991)]{cheng:robinson:1991:density}
B.~Cheng and P.M. Robinson.
\newblock Density estimation in strongly dependent non-linear time series.
\newblock \emph{Statistica Sinica}, 1\penalty0 (2):\penalty0 335--359, 1991.

\bibitem[Cheng and Ho(2008)]{cheng:Ho:2008:berry}
T.~Cheng and H.~Ho.
\newblock On {B}erry--{E}sseen bounds for non-instantaneous filters of linear
  processes.
\newblock \emph{Bernoulli}, 14\penalty0 (2):\penalty0 301--321, 2008.

\bibitem[Clausel et~al.(2012)Clausel, Roueff, Taqqu, and
  Tudor]{clausel:2012:large}
M.~Clausel, F.~Roueff, M.S. Taqqu, and C.~Tudor.
\newblock Large scale behavior of wavelet coefficients of non-linear
  subordinated processes with long memory.
\newblock \emph{Applied and Computational Harmonic Analysis}, 32\penalty0
  (2):\penalty0 223--241, 2012.

\bibitem[Clausel et~al.(2014)Clausel, Roueff, Taqqu, and
  Tudor]{clausel:2014:wavelet}
M.~Clausel, F.~Roueff, M.S. Taqqu, and C.~Tudor.
\newblock Wavelet estimation of the long memory parameter for hermite
  polynomial of {G}aussian processes.
\newblock \emph{ESAIM: Probability and Statistics}, 18:\penalty0 42--76, 2014.

\bibitem[Cs{\"o}rg{\"{o}}(2002)]{csorgo:2002:smoothing}
S.~Cs{\"o}rg{\"{o}}.
\newblock The smoothing dichotomy in nonparametric regression under long-memory
  errors.
\newblock \emph{Statistica neerlandica}, 56\penalty0 (2):\penalty0 132--142,
  2002.

\bibitem[Cs{\"o}rg{\"o} and Mielniczuk(1995)]{csorgo:mielniczuk:1995:density}
S.~Cs{\"o}rg{\"o} and J.~Mielniczuk.
\newblock Density estimation under long-range dependence.
\newblock \emph{The Annals of Statistics}, pages 990--999, 1995.

\bibitem[Cs{\"o}rg{\"o} and Mielniczuk(1999)]{csorgo:mielniczuk:1999:random}
S.~Cs{\"o}rg{\"o} and J.~Mielniczuk.
\newblock Random-design regression under long-range dependent errors.
\newblock \emph{Bernoulli}, pages 209--224, 1999.

\bibitem[Dehling and Philipp(2002)]{dehling:philipp:2002:empirical}
H.~Dehling and W.~Philipp.
\newblock Empirical process techniques for dependent data.
\newblock In \emph{Empirical process techniques for dependent data}, pages
  3--113. Springer, 2002.

\bibitem[Dehling and Taqqu(1989)]{dehling:taqqu:1989:empirical}
H.~Dehling and M.S. Taqqu.
\newblock The empirical process of some long-range dependent sequences with an
  application to {U}-statistics.
\newblock \emph{The Annals of Statistics}, pages 1767--1783, 1989.

\bibitem[Dehling and Taqqu(1991)]{dehling:taqqu:1991:bivariate}
H.~Dehling and M.S. Taqqu.
\newblock Bivariate symmetric statistics of long-range dependent observations.
\newblock \emph{Journal of Statistical Planning and Inference}, 28\penalty0
  (2):\penalty0 153--165, 1991.

\bibitem[Dehling et~al.(2002)Dehling, Mikosch, and
  Sorensen]{dehling:mikosch:2002:empirical}
H.~Dehling, T.~Mikosch, and M.~(editors) Sorensen.
\newblock \emph{Empirical process techniques for dependent data}.
\newblock Springer, 2002.

\bibitem[Dehling et~al.(2013)Dehling, Rooch, and Taqqu]{dehling:rooch:2013:non}
H.~Dehling, A.~Rooch, and M.S. Taqqu.
\newblock Non-parametric change-point tests for long-range dependent data.
\newblock \emph{Scandinavian Journal of Statistics}, 40\penalty0 (1):\penalty0
  153--173, 2013.

\bibitem[Denaranjo(1993)]{denaranjo:1993:non}
M.V.S. Denaranjo.
\newblock Non-central limit theorems for non-linear functionals of k {G}aussian
  fields.
\newblock \emph{Journal of multivariate analysis}, 44\penalty0 (2):\penalty0
  227--255, 1993.

\bibitem[Dobrushin and Major(1979)]{dobrushin:major:1979:non}
R.L. Dobrushin and P.~Major.
\newblock Non-central limit theorems for non-linear functional of {G}aussian
  fields.
\newblock \emph{Probability Theory and Related Fields}, 50\penalty0
  (1):\penalty0 27--52, 1979.

\bibitem[Doukhan(2003)]{doukhan:2003:models}
P.~Doukhan.
\newblock Models, inequalities, and limit theorems for stationary sequences.
\newblock In \emph{Theory and Applications of Long-Range Dependence}, pages
  43--100. Birkh{\"a}user, 2003.

\bibitem[Doukhan and Louhichi(1999)]{doukhan:louhichi:1999:new}
P.~Doukhan and S.~Louhichi.
\newblock A new weak dependence condition and applications to moment
  inequalities.
\newblock \emph{Stochastic Processes and their Applications}, 84\penalty0
  (2):\penalty0 313--342, 1999.

\bibitem[Fa{\"y} et~al.(2009)Fa{\"y}, Moulines, Roueff, and
  Taqqu]{fay:2009:estimators}
G.~Fa{\"y}, E.~Moulines, F.~Roueff, and M.S. Taqqu.
\newblock Estimators of long-memory: Fourier versus wavelets.
\newblock \emph{Journal of econometrics}, 151\penalty0 (2):\penalty0 159--177,
  2009.

\bibitem[Fox and Taqqu(1986)]{fox:taqqu:1986:large}
R.~Fox and M.S. Taqqu.
\newblock Large-sample properties of parameter estimates for strongly dependent
  stationary {G}aussian time series.
\newblock \emph{The Annals of Statistics}, pages 517--532, 1986.

\bibitem[Giraitis and Surgailis(1990)]{giraitis:surgailis:1990:central}
L.~Giraitis and D.~Surgailis.
\newblock A central limit theorem for quadratic forms in strongly dependent
  linear variables and its application to asymptotical normality of {W}hittle's
  estimate.
\newblock \emph{Probability Theory and Related Fields}, 86\penalty0
  (1):\penalty0 87--104, 1990.

\bibitem[Giraitis and Taqqu(1999)]{giraitis:taqqu:1999:whittle}
L.~Giraitis and M.S. Taqqu.
\newblock Whittle estimator for finite-variance non-{G}aussian time series with
  long memory.
\newblock \emph{Annals of Statistics}, pages 178--203, 1999.

\bibitem[Giraitis et~al.(2004)Giraitis, Leipus, Robinson, and
  Surgailis]{giraitis:2004:larch}
L.~Giraitis, R.~Leipus, P.M. Robinson, and D.~Surgailis.
\newblock {LARCH}, leverage, and long memory.
\newblock \emph{Journal of Financial Econometrics}, 2\penalty0 (2):\penalty0
  177--210, 2004.

\bibitem[Giraitis et~al.(2012)Giraitis, Koul, and
  Surgailis]{giraitis:koul:surgailis:2009:large}
L.~Giraitis, H.L. Koul, and D.~Surgailis.
\newblock \emph{Large Sample Inference for Long Memory Processes}.
\newblock World Scientific Publishing Company Incorporated, 2012.

\bibitem[Granger and Joyeux(1980)]{granger:joyeux:1980introduction}
C.W.J. Granger and R.~Joyeux.
\newblock An introduction to long-memory time series models and fractional
  differencing.
\newblock \emph{Journal of time series analysis}, 1\penalty0 (1):\penalty0
  15--29, 1980.

\bibitem[Guo and Koul(2007)]{guo:koul:2007:nonparametric}
H.~Guo and H.L. Koul.
\newblock Nonparametric regression with heteroscedastic long memory errors.
\newblock \emph{Journal of Statistical Planning and Inference}, 137\penalty0
  (2):\penalty0 379--404, 2007.

\bibitem[Hall et~al.(1998)Hall, Jing, and Lahiri]{hall:1998:sampling}
P.~Hall, B-Y Jing, and S.N. Lahiri.
\newblock On the sampling window method for long-range dependent data.
\newblock \emph{Statistica Sinica}, 8\penalty0 (4):\penalty0 1189--1204, 1998.

\bibitem[Herrndorf(1984)]{herrndorf:1984:functional}
N.~Herrndorf.
\newblock A functional central limit theorem for weakly dependent sequences of
  random variables.
\newblock \emph{The Annals of Probability}, pages 141--153, 1984.

\bibitem[Hidalgo(1997)]{hidalgo:1997:non}
J.~Hidalgo.
\newblock Non-parametric estimation with strongly dependent multivariate time
  series.
\newblock \emph{Journal of Time Series Analysis}, 18\penalty0 (2):\penalty0
  95--122, 1997.

\bibitem[Ho(1996)]{ho:1996:central}
H.~Ho.
\newblock On central and non-central limit theorems in density estimation for
  sequences of long-range dependence.
\newblock \emph{Stochastic processes and their Applications}, 63\penalty0
  (2):\penalty0 153--174, 1996.

\bibitem[Ho and Hsing(1997)]{ho:hsing:1997:limit}
H.~Ho and T.~Hsing.
\newblock Limit theorems for functionals of moving averages.
\newblock \emph{The Annals of Probability}, 25\penalty0 (4):\penalty0
  1636--1669, 1997.

\bibitem[Ho and Sun(1987)]{ho:sun:1987:central}
H.~Ho and T.~Sun.
\newblock A central limit theorem for non-instantaneous filters of a stationary
  {G}aussian process.
\newblock \emph{Journal of Multivariate Analysis}, 22\penalty0 (1):\penalty0
  144--155, 1987.

\bibitem[Hosking(1996)]{hosking:1996:asymptotic}
J.R.M. Hosking.
\newblock Asymptotic distributions of the sample mean, autocovariances, and
  autocorrelations of long-memory time series.
\newblock \emph{Journal of Econometrics}, 73\penalty0 (1):\penalty0 261--284,
  1996.

\bibitem[Ibragimov(1962)]{ibragimov:1962:some}
I.A. Ibragimov.
\newblock Some limit theorems for stationary processes.
\newblock \emph{Theory of Probability \& Its Applications}, 7\penalty0
  (4):\penalty0 349--382, 1962.

\bibitem[Jach et~al.(2012)Jach, McElroy, and Politis]{jach:2012:subsampling}
A.~Jach, T.~McElroy, and D.N. Politis.
\newblock Subsampling inference for the mean of heavy-tailed long-memory time
  series.
\newblock \emph{Journal of Time Series Analysis}, 33\penalty0 (1):\penalty0
  96--111, 2012.

\bibitem[Janson(1997)]{janson:1997:gaussian}
S.~Janson.
\newblock \emph{Gaussian Hilbert Spaces}.
\newblock Cambridge Tracts in Mathematics. Cambridge University Press, 1997.

\bibitem[Koul and Stute(1998)]{koul:1998:regression}
H.L. Koul and W.~Stute.
\newblock Regression model fitting with long memory errors.
\newblock \emph{Journal of statistical planning and inference}, 71\penalty0
  (1):\penalty0 35--56, 1998.

\bibitem[K\"unsch(1987)]{kunsch:1987:statistical}
H.R. K\"unsch.
\newblock Statistical aspects of self-similar processes.
\newblock In \emph{Proceedings of the first World Congress of the Bernoulli
  Society}, volume~1, pages 67--74. VNU Science Press Utrecht, 1987.

\bibitem[L{\'e}vy-Leduc and Taqqu(2013)]{levy-leduc:taqqu:2013:long}
C.~L{\'e}vy-Leduc and M.S. Taqqu.
\newblock Long-range dependence and the rank of decompositions.
\newblock \emph{Fractal Geometry and Dynamical Systems in Pure and Applied
  Mathematics II: Fractals in Applied Mathematics}, 601:\penalty0 289, 2013.

\bibitem[L{\'e}vy-Leduc et~al.(2011)L{\'e}vy-Leduc, Boistard, Moulines, Taqqu,
  and Reisen]{levy:boistard:taqqu:reisen:2011:asymptotic}
C.~L{\'e}vy-Leduc, H.~Boistard, E.~Moulines, M.S. Taqqu, and V.A. Reisen.
\newblock Asymptotic properties of {U}-processes under long-range dependence.
\newblock \emph{The Annals of Statistics}, 39\penalty0 (3):\penalty0
  1399--1426, 2011.

\bibitem[L\'evy-Leduc et~al.(2011)L\'evy-Leduc, Boistard, Moulines, and
  Taqqu]{levy-leduc:2010:robust}
C~L\'evy-Leduc, H.~Boistard, E.~Moulines, and Reisen~V.A. Taqqu, M.S.
\newblock Robust estimation of the scale and of the autocovariance function of
  {G}aussian short and long-range dependent processes.
\newblock \emph{Journal of Time Series Analysis}, 32:\penalty0 135--156, 2011.

\bibitem[Major(2014)]{major:2014:multiple}
P.~Major.
\newblock \emph{Multiple Wiener-It{\^o} Integrals: With Applications to Limit
  Theorems}.
\newblock Lecture Notes in Mathematics. Springer, 2nd edition, 2014.

\bibitem[Mandelbrot and Wallis(1969)]{mandelbrot:wallis:1969:some}
B.B. Mandelbrot and J.R. Wallis.
\newblock Some long-run properties of geophysical records.
\newblock \emph{Water resources research}, 5\penalty0 (2):\penalty0 321--340,
  1969.

\bibitem[Masry and Mielniczuk(1999)]{masry:mielniczuk:1999:local}
E.~Masry and J.~Mielniczuk.
\newblock Local linear regression estimation for time series with long-range
  dependence.
\newblock \emph{Stochastic Processes and their Applications}, 82\penalty0
  (2):\penalty0 173--193, 1999.

\bibitem[Nordman and Lahiri(2005)]{nordman:lahiri:2005:validity}
D.J. Nordman and S.N. Lahiri.
\newblock Validity of the sampling window method for long-range dependent
  linear processes.
\newblock \emph{Econometric Theory}, 21\penalty0 (06):\penalty0 1087--1111,
  2005.

\bibitem[Pelletier and Turcotte(1997)]{pelletier:turcotte:1997:long}
J.D. Pelletier and D.L. Turcotte.
\newblock Long-range persistence in climatological and hydrological time
  series: analysis, modeling and application to drought hazard assessment.
\newblock \emph{Journal of Hydrology}, 203\penalty0 (1):\penalty0 198--208,
  1997.

\bibitem[Phillips(1987)]{phillips:1987:towards}
P.C.B. Phillips.
\newblock Towards a unified asymptotic theory for autoregression.
\newblock \emph{Biometrika}, pages 535--547, 1987.

\bibitem[Pipiras and Taqqu(2010)]{pipiras:taqqu:2010:regularization}
V.~Pipiras and M.S. Taqqu.
\newblock Regularization and integral representations of {H}ermite processes.
\newblock \emph{Statistics and Probability Letters}, 80\penalty0 (23):\penalty0
  2014--2023, 2010.

\bibitem[Pipiras and Taqqu(2017)]{pipiras:taqqu:2017:long}
V.~Pipiras and M.S. Taqqu.
\newblock \emph{Long-Range Dependence and Self-Similarity}.
\newblock Cambridge University Press, 2017.

\bibitem[Psaradakis(2010)]{psaradakis:2010:inference}
Z.~Psaradakis.
\newblock On inference based on the one-sample sign statistic for long-range
  dependent data.
\newblock \emph{Computational Statistics}, 25\penalty0 (2):\penalty0 329--340,
  2010.

\bibitem[Robinson(1995)]{robinson:1995:gaussian}
P.M. Robinson.
\newblock {G}aussian semiparametric estimation of long range dependence.
\newblock \emph{The Annals of Statistics}, pages 1630--1661, 1995.

\bibitem[Rosenblatt(1956)]{rosenblatt:1956:central}
M.~Rosenblatt.
\newblock A central limit theorem and a strong mixing condition.
\newblock \emph{Proceedings of the National Academy of Sciences of the United
  States of America}, 42\penalty0 (1):\penalty0 43, 1956.

\bibitem[Rosenblatt(1961)]{rosenblatt:1961:independence}
M.~Rosenblatt.
\newblock Independence and dependence.
\newblock In \emph{Proc. Fourth Berkeley Symp. Math. Statist. Probab},
  volume~2, pages 431--443, 1961.

\bibitem[Samorodnitsky(2016)]{samorodnitsky:2016:stochastic}
G.~Samorodnitsky.
\newblock \emph{Stochastic Processes and Long Range Dependence}.
\newblock Springer, 2016.

\bibitem[Schweingruber(1996)]{schweingruber:1996:tree}
F.H. Schweingruber.
\newblock \emph{Tree rings and environment: dendroecology.}
\newblock Paul Haupt AG Bern, 1996.

\bibitem[Shao(2010)]{shao:2010:self}
X.~Shao.
\newblock A self-normalized approach to confidence interval construction in
  time series.
\newblock \emph{Journal of the Royal Statistical Society: Series B (Statistical
  Methodology)}, 72\penalty0 (3):\penalty0 343--366, 2010.

\bibitem[Shao(2011)]{shao:2011:simple}
X.~Shao.
\newblock A simple test of changes in mean in the possible presence of
  long-range dependence.
\newblock \emph{Journal of Time Series Analysis}, 32\penalty0 (6):\penalty0
  598--606, 2011.

\bibitem[Surgailis(1982)]{surgailis:1982:zones}
D.~Surgailis.
\newblock Zones of attraction of self-similar multiple integrals.
\newblock \emph{Lithuanian Mathematical Journal}, 22\penalty0 (3):\penalty0
  327--340, 1982.

\bibitem[Surgailis(2000)]{surgailis:2000:long}
D.~Surgailis.
\newblock Long-range dependence and {A}ppell rank.
\newblock \emph{Annals of Probability}, pages 478--497, 2000.

\bibitem[Taqqu(1975)]{taqqu:1975:weak}
M.S. Taqqu.
\newblock Weak convergence to fractional {B}rownian motion and to the
  {R}osenblatt process.
\newblock \emph{Probability Theory and Related Fields}, 31\penalty0
  (4):\penalty0 287--302, 1975.

\bibitem[Taqqu(1979)]{taqqu:1979:convergence}
M.S. Taqqu.
\newblock Convergence of integrated processes of arbitrary {H}ermite rank.
\newblock \emph{Probability Theory and Related Fields}, 50\penalty0
  (1):\penalty0 53--83, 1979.

\bibitem[Taqqu and Teverovsky(1997)]{taqqu:teverovsky:1997:robustness}
M.S. Taqqu and V.~Teverovsky.
\newblock Robustness of {W}hittle-type estimators for time series with
  long-range dependence.
\newblock \emph{Communications in Statistics. Stochastic Models}, 13\penalty0
  (4):\penalty0 723--757, 1997.

\bibitem[Taqqu et~al.(1995)Taqqu, Teverovsky, and
  Willinger]{taqqu:1995:estimators}
M.S. Taqqu, V.~Teverovsky, and W.~Willinger.
\newblock Estimators for long-range dependence: an empirical study.
\newblock \emph{Fractals}, 3\penalty0 (04):\penalty0 785--798, 1995.

\bibitem[Terrin and Taqqu(1990)]{terrin:taqqu:1990:noncentral}
N.~Terrin and M.S. Taqqu.
\newblock A noncentral limit theorem for quadratic forms of {G}aussian
  stationary sequences.
\newblock \emph{Journal of Theoretical Probability}, 3\penalty0 (3):\penalty0
  449--475, 1990.

\bibitem[Wu(2002)]{wu:2002:central}
W.B. Wu.
\newblock Central limit theorems for functionals of linear processes and their
  applications.
\newblock \emph{Statistica Sinica}, 12\penalty0 (2):\penalty0 635--650, 2002.

\bibitem[Wu(2005)]{wu:2005:nonlinear}
W.B. Wu.
\newblock Nonlinear system theory: Another look at dependence.
\newblock \emph{Proceedings of the National Academy of Sciences of the United
  States of America}, 102\penalty0 (40):\penalty0 14150--14154, 2005.

\bibitem[Wu(2006)]{wu:2006:unit}
W.B. Wu.
\newblock Unit root testing for functionals of linear processes.
\newblock \emph{Econometric Theory}, 22\penalty0 (01):\penalty0 1--14, 2006.

\bibitem[Wu and Mielniczuk(2002)]{wu:mielniczuk:2002:kernel}
W.B. Wu and J.~Mielniczuk.
\newblock Kernel density estimation for linear processes.
\newblock \emph{Annals of Statistics}, pages 1441--1459, 2002.

\bibitem[Wu et~al.(2010)Wu, Huang, and Zheng]{wu:huang:2010:covariance}
W.B. Wu, Y.~Huang, and W.~Zheng.
\newblock Covariances estimation for long-memory processes.
\newblock \emph{Advances in Applied Probability}, 42\penalty0 (1):\penalty0
  137--157, 2010.

\bibitem[Zhang et~al.(2013)Zhang, Ho, Wendler, and Wu]{zhang:2013:block}
T.~Zhang, H-C Ho, M.~Wendler, and W.B. Wu.
\newblock Block sampling under strong dependence.
\newblock \emph{Stochastic Processes and their Applications}, 123\penalty0
  (6):\penalty0 2323--2339, 2013.

\bibitem[Zhao et~al.(2010)Zhao, Tian, and Xia]{zhao:2010:ratio}
W.~Zhao, Z.~Tian, and Z.~Xia.
\newblock Ratio test for variance change point in linear process with long
  memory.
\newblock \emph{Statistical Papers}, 51\penalty0 (2):\penalty0 397--407, 2010.

\end{thebibliography}

\end{document}